\documentclass[review,hidelinks,onefignum,onetabnum]{siamart220329}

\usepackage{lipsum}
\usepackage{amsfonts}
\usepackage{graphicx}
\usepackage{epstopdf}

\nolinenumbers

\newtheorem{example}{Example}

\usepackage{enumerate}
\usepackage{comment}
\usepackage{csquotes}
\usepackage{subcaption}
\usepackage{algorithm}
\usepackage{algorithmicx} 
\usepackage{algpseudocode}

\graphicspath{{figure_eps/}}

\ifpdf
  \DeclareGraphicsExtensions{.eps,.pdf,.png,.jpg}
\else
  \DeclareGraphicsExtensions{.eps}
\fi

\newsiamremark{remark}{Remark}
\newsiamremark{hypothesis}{Hypothesis}
\crefname{hypothesis}{Hypothesis}{Hypotheses}
\newsiamthm{claim}{Claim}

\headers{Machine learning optimization of M\"untz spectral methods}{Wei Zeng, Chuanju Xu, Yiming Lu and Qian Wang}

\title{Machine learning-based parameter optimization for M\"untz spectral methods 
\thanks{Submitted to the editors DATE.
\funding{This work was supported by the National Natural Science Foundation of China (Grants 12372284, 12371408 and U2230402) and the Postdoctoral Fellowship Program of CPSF (Grant GZC20230213).}}}
\author{Wei Zeng\thanks{Mechanics Division, Beijing Computational Science Research Center, 100193 Beijing, China
  (\email{zengwei@csrc.ac.cn}, \email{qian.wang@csrc.ac.cn}).}
\and Chuanju Xu\thanks{School of Mathematical Sciences and
	Fujian Provincial Key Laboratory of Mathematical Modeling and High Performance Scientific Computing, Xiamen University, 361005 Xiamen, China
  (\email{cjxu@xmu.edu.cn}).}
\and Yiming Lu\thanks{Institute for Aero Engine, Tsinghua University, 100084 Beijing, China
	(\email{lu-ym24@mails.tsinghua.edu.cn}).}
\and Qian Wang\footnotemark[2]}

\usepackage{amsopn}

\ifpdf
\hypersetup{
  pdftitle={Machine learning-based parameter optimization for M\"untz spectral methods},
  pdfauthor={Wei Zeng, Chuanju Xu, Yiming Lu and Qian Wang}
}
\fi

\begin{document}

\maketitle
\begin{abstract}
Spectral methods employing non-standard polynomial bases, such as M\"untz polynomials, have proven effective for accurately solving problems with solutions exhibiting low regularity, notably including sub-diffusion equations. However, due to the absence of theoretical guidance, the key parameters controlling the exponents of M\"untz polynomials are usually determined empirically through extensive numerical experiments, leading to a time-consuming tuning process.
To address this issue, we propose a novel machine learning-based optimization framework for the M\"untz spectral method. As an illustrative example, we optimize the parameter selection for solving time-fractional partial differential equations (PDEs). Specifically, an artificial neural network (ANN) is employed to predict optimal parameter values based solely on the time-fractional order as input. The ANN is trained by minimizing solution errors on a one-dimensional time-fractional convection-diffusion equation featuring manufactured exact solutions that manifest singularities of varying intensity, covering a comprehensive range of sampled fractional orders.
Numerical results for time-fractional PDEs in both one and two dimensions demonstrate that the ANN-based parameter prediction significantly improves the accuracy of the M\"untz spectral method. 
Moreover, the trained ANN generalizes effectively from one-dimensional to two-dimensional cases, highlighting its robustness across spatial dimensions. 
Additionally, we verify that the ANN substantially outperforms traditional function approximators, such as spline interpolation, in both prediction accuracy and training efficiency.
The proposed optimization framework can be extended beyond fractional PDEs, offering a versatile and powerful approach for spectral methods applied to various low-regularity problems.

\end{abstract}



\section{Introduction}
\label{sec1}

Partial differential equations (PDEs) are extensively used to model complex problems across a wide range of scientific disciplines \cite{evans2022partial}. 
Spectral methods \cite{Shen2011Spectral} are a class of techniques developed to numerically solve PDEs, in which the solution is usually approximated by an expansion of orthogonal basis functions, such as trigonometric or polynomial functions. 
These basis functions are employed to accurately capture the key characteristics of solutions in areas like fluid dynamics, heat transfer, wave propagation, and quantum mechanics \cite{karniadakis2005spectral,gottlieb1977numerical}. The Fourier spectral method, which relies on trigonometric basis functions, is limited to problems with periodic boundary conditions. In contrast, spectral methods based on orthogonal polynomial basis functions provide greater flexibility and can be applied to a wider variety of boundary conditions and problem domains.


Classical spectral methods are renowned for their high accuracy, often achieving exponential convergence rates when applied to problems with sufficiently smooth solutions. However, these methods tend to perform poorly with non-smooth solutions, as smooth basis functions like polynomials and trigonometric functions are not well-suited to capturing abrupt changes or singularities. The M\"untz spectral method, 
which allows for the selection of fractional basis functions, introduces adaptability that can be tailored to the smoothness and specific features of the solution. This flexibility improves both efficiency and accuracy, particularly for problems involving non-smooth solutions. The M\"untz spectral method has been successfully applied to address challenges associated with singularities, such as the Poisson equation with mixed Dirichlet–Neumann boundary conditions \cite{Shen2016Muntz}, Volterra-type integral equations with weakly singular kernels \cite{Hou2019Muntz}, and time-fractional PDEs \cite{Hou2018Muntz,Hou2017A,Zeng2024An}. 
Here we are particularly interested in
sub-diffusion equations, which departs from the classical Fickian framework of Brownian motion. Anomalous diffusion is frequently observed across diverse physical and biological systems, particularly in contexts where particle trapping and binding occur \cite{BWZW,FBSBW,GW,SLSJ,SMWFC,SL,SM,MKW,AMYPL,Nig,Mai95}, making time-fractional PDE a valuable tool for capturing such complex behaviors.

In our M\"untz spectral method, we consider a simple case in which the approximation space is spanned by 
$\{x^{\lambda_0}, x^{\lambda_1}, \dots, x^{\lambda_N}\}$, 
$\lambda_n=n\lambda$ with $\lambda$ being a free parameter. It has been reported (see, e.g., \cite{Hou2018Muntz,Hou2019Muntz}) that the choice of $\lambda$ has a significant impact on the accuracy of the method. If the fractional polynomial basis functions are poorly aligned with the characteristics of the solution, the convergence rate may decline, leading to a substantial loss of accuracy \cite{Zeng2024An}. Therefore, it is crucial to set the parameter $\lambda$ appropriately to ensure accurate simulations. 
However, there is a lack of theoretical guidance for determining the optimal value of $\lambda$, as it is difficult to assess the solution’s regularity before simulation. Additionally, the relationship between the solution’s regularity and the optimal $\lambda$ is highly nonlinear and complex. As a result, $\lambda$ is usually chosen empirically by testing multiple candidate values and selecting the one that yields the highest accuracy \cite{Shen2016Muntz,Hou2018Muntz,Hou2017A,Hou2019Muntz}. This time-consuming parameter tuning process has become a major bottleneck in the practical implementation of the M\"untz spectral method.

In this work, a machine learning approach is developed to determine the optimal value of the free parameter $\lambda$, addressing the bottleneck of manual parameter selection and enhancing the accuracy of the M\"untz spectral method. 
Specifically, the M\"untz spectral method for time-fractional PDEs is optimized to demonstrate the proposed machine learning approach. An artificial neural network (ANN) is trained to predict the optimal parameter $\lambda$ for the M\"untz spectral method, with the time-fractional order $\mu$ of the PDE as input.
In the M\"untz spectral method, both the time-fractional order and the forcing function play key roles in determining the optimal value of $\lambda$. In this paper, we focus on training the ANN to approximate the mapping from the time-fractional order $\mu$ to the optimal value of $\lambda$, with the goal of applying the network to a broad class of time-fractional PDEs with varying forcing functions.

The ANN is trained by minimizing the solution errors of a one-dimensional time-fractional convection-diffusion equation. To facilitate this, a set of time-fractional orders is sampled from the parameter domain $\left(0,1\right)$, and a set of manufactured solutions with various singularities is generated to serve as exact solutions. This process yields a dataset of time-fractional convection-diffusion equations, based on the sampled time-fractional orders and manufactured solutions.
For each equation in the dataset, a numerical solution is computed using the M\"untz spectral method, where the parameter $\lambda$ is predicted by the ANN. The loss function is defined as the mean Frobenius norm of the errors between these numerical solutions and their corresponding exact solutions. 
During training, the ANN learns to predict the optimal parameter values for the given time-fractional PDEs by minimizing the loss function. 
Furthermore, in practical applications, the optimal parameter value can be pre-computed, ensuring no additional computational cost during runtime and thus preserving the efficiency of the M\"untz spectral method.
The trained ANN is integrated with the M\"untz spectral method to solve time-fractional PDEs in both one and two dimensions. Numerical results indicate that the ANN-driven parameter prediction substantially improves solution accuracy. Furthermore, although initially trained in a one-dimensional context, the network demonstrates strong generalization to two-dimensional time-fractional problems, underscoring its robust performance across different spatial dimensions. 

It should be noted that, in addition to artificial neural networks (ANN), various other trainable function approximators may serve as parameter prediction models within the proposed optimization framework. Alternative parameter prediction models can also be integrated into the M\"untz spectral solver and trained by minimizing solution errors associated with the one-dimensional time-fractional PDE. For comparison purposes, we also train a cubic spline interpolation model, a classical method frequently used in numerical approximation. Numerical results demonstrate that the ANN significantly outperforms spline interpolation, exhibiting superior prediction accuracy and training efficiency.

The key contributions and novel findings of this paper include:
\renewcommand{\labelenumi}{(\arabic{enumi})} 
\begin{enumerate}
	\item An optimization approach based on machine learning is proposed for the M\"untz spectral method. 
	Specifically, an ANN is employed to predict the optimal value of the free parameter $\lambda$, which determines the exponents of the M\"untz polynomial basis functions.
	As an illustrative example, we optimize the parameter selection for solving time-fractional PDEs, where the time-fractional order $\mu$ is taken as the input of the ANN. 
	
	\item The ANN is trained by minimizing the solution errors of a one-dimensional time-fractional convection-diffusion equation. To this end, a dataset is generated by sampling a range of time-fractional orders and constructing manufactured exact solutions with varying singularities.
	For each equation in the dataset, a numerical solution is computed using the M\"untz spectral method, with the parameter predicted by the ANN.  The ANN training involves minimizing a loss function defined as the mean Frobenius norm of the errors between these numerical solutions and their corresponding exact solutions.
	\item Numerical results demonstrate that the ANN-driven parameter prediction significantly enhances the accuracy of the M\"untz spectral method. Moreover, despite being initially trained on one-dimensional problems, the ANN exhibits strong generalization to two-dimensional time-fractional problems, highlighting its robust performance across spatial dimensions.
	
	\item The proposed machine learning optimization approach can be extended to a broader range of low-regularity problems beyond those involving fractional operators. In these contexts, the ANN adaptively accounts for the solution’s regularity, enabling more accurate simulations through the efficient prediction of optimal parameters for numerical schemes. 
\end{enumerate}

The remainder of this paper is organized as follows. Section \ref{sec2} introduces the M\"untz spectral method for time-fractional PDEs.
Section \ref{sec3} presents the machine learning optimization of the M\"untz spectral method. Numerical results are given in
Section \ref{sec4} and concluding remarks are given in Section \ref{sec5}.

\section{M\"untz spectral method}
\label{sec2}
\subsection{Fractional polynomial basis}
Let $I=(0,1)$. The classical Weierstrass theorem states that every continuous function on a closed interval $I$ can be uniformly approximated as closely as desired by a polynomial function. Later, this theorem was generalized by M\"untz and Sz\'asz, who proved that the space $\rm{span}$$\{x^{\lambda_0}, x^{\lambda_1},x^{\lambda_2}, \cdots\}$, called M\"untz polynomial space, is dense in $L^2(I)$ if and only if $\sum^\infty_{n=0}\lambda_n^{-1}=+\infty$, where $\left\{0= \lambda_0<\lambda_1<\lambda_2<\cdots\right\}$ is an increasing sequence of real numbers. This generalization, usually called the M\"untz-Sz\'asz theorem \cite{Ortiz2005Herman,Almira2007Muntz}, is significant since $L^2$-spaces include functions that are not necessarily continuous but are square-integrable over the interval. According to this theorem, the traditional polynomial space $\mathrm{span}$$\left\{1, x,x^2,\cdots\right\}$ can be extended to the fractional polynomial space $\mathrm{span}$$\left\{1, x^{\lambda_1},x^{\lambda_2},\cdots \right\}$, for function approximation.
This class of fractional polynomials serves as a foundational basis in approximation theory and underpins various numerical methods, particularly spectral methods \cite{Shen2016Muntz, Hou2018Muntz,Hou2017A, Cui2024Muntz}. 
Spectral methods that utilize M\"untz polynomial bases are known as M\"untz spectral methods.


The fractional polynomial basis $\{1,x^{\lambda_1},x^{\lambda_2}, \cdots\}$, with varying exponents $\{0, \lambda_1,\lambda_2, \cdots\}$, can be employed to approximate non-smooth solutions to specific problems \cite{Cui2024Muntz}. In this paper, we consider a particular class of fractional polynomial basis with $\lambda_n = n\lambda$, where $\lambda \in (0,1]$ is a  positive number. The corresponding M\"untz polynomials are capable of capturing left endpoint singularities, which is a typical behavior of solutions to sub-diffusion equations.

\begin{figure}[htbp!]
	\centering
	\begin{subfigure}{0.48\textwidth}
		\includegraphics[width=\linewidth]{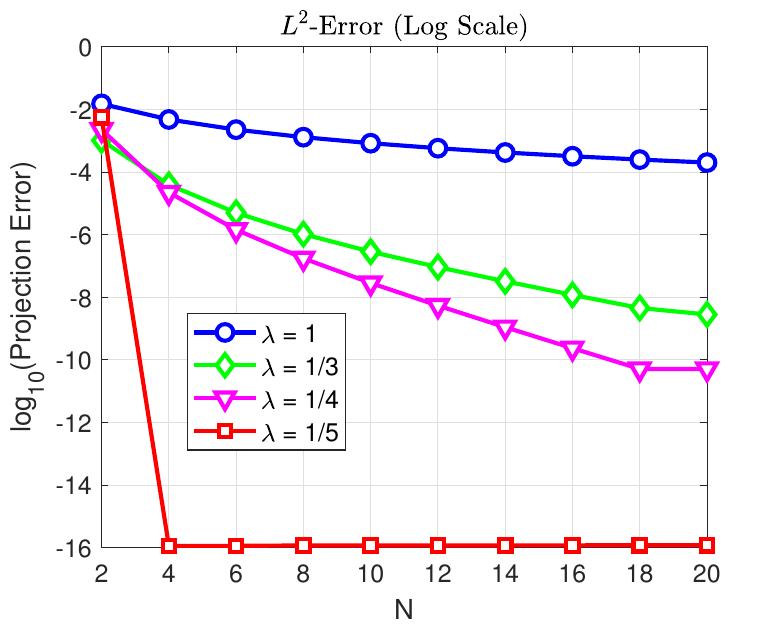}
		\caption{$f=x^{3/5}$}
		\label{pro_1}
	\end{subfigure}
	\begin{subfigure}{0.48\textwidth}
		\includegraphics[width=\linewidth]{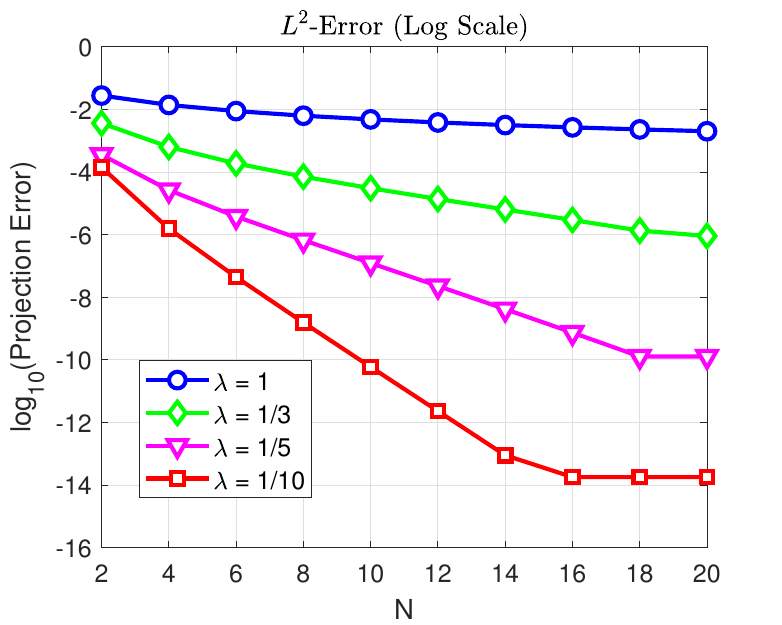}
		\caption{$f=x^{\sqrt{3}/10}$}
		\label{pro_2}
	\end{subfigure}
	\caption{Projection error versus the polynomial degree $N$.}
	\label{fig:pro_error}
\end{figure}

Here, we present an experimental investigation into the impact of the parameter $\lambda$ on the accuracy of approximating non-smooth functions using fractional polynomials. Our numerical experiments focus on two representative test functions, $f\left(x\right)= x^{3/5}$ and $f\left(x\right)= x^{\sqrt{3}/10}$. Fig. \ref{fig:pro_error} offers a detailed comparison of projection errors across various $\lambda$ values.
The results reveal that traditional polynomial approximations (where $\lambda=1$) display relatively slow convergence. In contrast, the M\"untz polynomial approximation with $\lambda_n= n \lambda$ achieves significantly higher accuracy and faster convergence rates, demonstrating its enhanced efficacy for non-smooth functions. Specifically, Fig. \ref{pro_1} shows optimal convergence behavior at $\lambda=1/5$, attributed to the smooth mapping function $f(x^{1/\lambda})=x^3$.
However, for $f\left(x\right)= x^{\sqrt{3}/10}$, the irrational exponent $\sqrt{3}/10$ prevents the mapping function $f(x^{1/\lambda})$ from being smooth. Additionally, the regularity of the mapping function $f(x^{1/\lambda})$ can be refined by decreasing $\lambda$. Consequently, the optimal convergence behavior is observed for the smallest $\lambda$ among the four options presented in Fig. \ref{pro_2}.

The orthogonal M\"{u}ntz Jacobi polynomials are used in this work to efficiently implement the spectral methods. Following \cite{Hou2018Muntz,Hou2017A}, the M\"untz Jacobi polynomial basis, constructed by applying Gram-Schmidt orthogonalization to the fractional polynomial basis $\left\{1, x^\lambda, \cdots, x^{n\lambda}, \cdots \right\}$, is defined as follows.

\begin{definition}
	For $\alpha,\beta\geq -1,$ $0<\lambda\leq 1,$ the M\"untz Jacobi polynomial of degree $n+l$ on $I$ is 
	\begin{equation} \label{muntz-jacobi}
		J^{\alpha,\beta,\lambda}_{n+l}(t)=
		\begin{cases}
			J^{\alpha,\beta}_{n}(2t^{\lambda}-1),\ &\alpha,\beta>-1,\\
			\frac{n+\alpha+1}{n+1}t^{\lambda}J^{\alpha,1}_n(2t^{\lambda}-1),\ &\alpha>-1,\beta=-1,\\
			\frac{n+\beta+1}{n+1}(1-t^{\lambda})J^{1,\beta}_n(2t^{\lambda}-1),\ &\alpha=-1,\beta>-1,\\
			-(1-t^{\lambda})t^{\lambda}J^{1,1}_n(2t^{\lambda}-1),\ &\alpha=\beta=-1,
		\end{cases}
	\end{equation}
	where $J^{\alpha,\beta}_n(s)$ is the usual Jacobi polynomial for index $\{\alpha,\beta\}$ of degree $n$, and
	\begin{equation} \label{l}
		l=
		\begin{cases}
			0,\ &\alpha,\beta>-1,\\
			1,\ &\alpha>-1,\beta=-1~~\mathrm{or}~~\alpha=-1,\beta>-1,\\
			2,\ &\alpha=\beta=-1.
		\end{cases}
	\end{equation}
\end{definition}
When $\lambda=1$, the M\"untz Jacobi polynomials $J^{\alpha,\beta,1}_{n+l}(t)$ are reduced to the shifted generalized Jacobi polynomials that are orthogonal on $I$ with respect to the weight $(1-t)^{\alpha}t^{\beta}$. More generally, the M\"untz Jacobi polynomials $J^{\alpha,\beta,\lambda}_{n+l}(t)$ are mutually orthogonal on $I$ with respect to the weight $\chi^{\alpha,\beta,\lambda}(t)=\lambda(1-t^{\lambda})^{\alpha}t^{(\beta+1){\lambda}-1},$ i.e.,
\begin{equation}\label{OrthoP}
	\int^1_{0}J^{\alpha,\beta,\lambda}_{n+l}(t)J^{\alpha,\beta,\lambda}_{m+l}(t)\chi^{\alpha,\beta,\lambda}(t)dt=\gamma^{\alpha,\beta}_{n+l}\delta_{m+l,n+l},
\end{equation}
where $\delta_{m+l,n+l}$ is the Kronecker function, and 
\begin{equation*}
	\gamma^{\alpha,\beta}_{n+l}=
	\frac{\Gamma(n+l+\alpha+1)\Gamma(n+l+\beta+1)}{(2(n+l)+\alpha+\beta+1)(n+l)!\Gamma(n+l+\alpha+\beta+1)}.
\end{equation*}
Note that we use the M\"untz Jacobi polynomial basis with $\alpha=0.5$ and $\beta=-1$, to effectively approximate solutions that exhibit singularities in this paper.

\subsection{M\"untz spectral method for time-fractional PDEs}

Let $\Omega=(-1,1), \Omega^d:= (-1,1)^d$, where $d$ represents the dimensionality. 
We consider a time-fractional convection-diffusion equation
\begin{equation}\label{TFDE}
	D_t^\mu u(\boldsymbol x,t)=\kappa\Delta u(\boldsymbol x,t)-\rho \nabla  u(\boldsymbol x,t)+f(\boldsymbol x,t),\ \boldsymbol x\in \Omega^d,~t\in I,
\end{equation}
subject to a homogeneous Dirichlet boundary condition
\begin{equation}\label{Bou}
	u(\boldsymbol x,t)\big|_{\partial\Omega^d}=0,\ t\in I,
\end{equation}
and a homogeneous initial condition
\begin{equation}\label{Ini}
	u(\boldsymbol x,0)=0,\ \boldsymbol x\in \Omega^d.
\end{equation}
Here, $\mu \in (0,1)$ is the fractional order, $\kappa>0$ and $\rho>0$ are the diffusion and convection coefficients, respectively, and $f$ is a forcing function. $D^\mu_t$ is the Caputo fractional derivative defined by
\begin{equation*}
	D^\mu_t u=I^{1-\mu}_t(\partial_tu), 
\end{equation*}
where
\begin{equation*}
	I^{\mu}_tu(t)=\omega_{\mu}\ast u(t)=\int^t_0\omega_{\mu}(t-s)u(s)ds, \quad \omega_{\mu}(t)=
	\frac{t^{\mu-1}}{\Gamma(\mu)}.
\end{equation*}
In this work, we focus on problems in one dimension ($d=1$) and two dimensions ($d=2$). Specifically, for $d=1$, we set $\boldsymbol x :=x$, and for $d=2$, $\boldsymbol x :=(x,y)$. Without loss of generality, we only consider homogeneous initial conditions in this paper. Non-homogeneous initial conditions can be handled by standard homogenization. 
A space-time spectral method will be proposed to solve equation \eqref{TFDE}. Notably, the solution to \eqref{TFDE} is smooth in the spatial dimension but exhibits singularities near the initial time. Such low regularity can substantially reduce the accuracy of conventional spectral methods. 
Consequently, standard orthogonal basis functions, such as Legendre polynomials, are applied in the spatial domain, while M\"untz Jacobi polynomials are employed as basis functions in the time direction to effectively handle solution singularities near the initial time.

In the time direction, a set of $N+1$ M\"untz Jacobi polynomials 
\begin{equation*}
	J^{\alpha,-1,\lambda}_{n+1}(t), \quad \alpha>-1,\quad 0<\lambda \leq 1, \quad n=0,\cdots,N, 
\end{equation*}
are used as basis functions, to approximate the solution in space
\begin{equation*}
	S^{\alpha,-1}_{N,\lambda}(I):=\mathrm{span}\left\{J^{\alpha,-1,\lambda}_{n+1}(t)\big|n=0,\cdots,N \right\}=\left\{ v\in P^{\lambda}_{N}(I) \big|v(0)=0 \right\},
\end{equation*}
where 
\begin{equation*}
	P^{\lambda}_{N}(I):=\mathrm{span} \left\{1,t^{\lambda},t^{2\lambda},\cdots,t^{N\lambda}\right\}
\end{equation*}
is the fractional polynomial space.
In each dimension of the spatial domain, a set of $M-1$ algebraic polynomials
\begin{equation}\label{basis}
	\phi_m(x)=c_m \left[L_m(x)-L_{m+2}(x)\right], \quad c_m=\frac{1}{\sqrt{4m+6}}, \quad m=0,\cdots,M-2,
\end{equation}
are used as basis functions, where $L_n(x)$ is the Legendre polynomial of degree $n$, to approximate the solution in space
\begin{equation*}
	\mathcal{P}_M(\Omega)=\mathrm{span}\left\{\phi_0(x),\ \phi_1(x),\ \ldots,\ \phi_{M-2}(x)\right\}.
\end{equation*}
The inner products
\begin{equation*}
	s^x_{jk}=(\phi^{\prime}_k(x),\phi^{\prime}_j(x)),\ \ 
	m^x_{jk}=(\phi_k(x),\phi_j(x)),\ \ c^x_{jk}=(\partial_x\phi_k(x),\phi_j(x)),
\end{equation*}
can be directly computed \cite{Shen1994Efficient} as 
\begin{align}\label{ab}
	\begin{split}
		s^x_{jk}=
		\begin{cases}
			1,\ &k=j,\\
			0,\ &k\neq j,
		\end{cases}
		\ \ \ \
		&m^x_{jk}=m^x_{kj}=
		\begin{cases}
			c_kc_j(\frac{2}{2j+1}+\frac{2}{2j+5}),\ &k=j,\\
			-c_kc_j\frac{2}{2k+1},\ &k=j+2,\\
			0,\ &\mathrm{otherwise},
		\end{cases}\\
		&c^x_{jk}=-c^x_{kj}=
		\begin{cases}
			2c_kc_j, &k=j+1,\\
			0, &\mathrm{otherwise}.
		\end{cases}
	\end{split}
\end{align}
The above equation shows that $(s^x_{jk})$ is a diagonal matrix, $(m^x_{jk})$ is a pentadiagonal symmetric matrix, whereas $(c^x_{jk})$ is an antisymmetric matrix.

Let $Q=\Omega^d \times I$. We define a basis size tuple $L$ and a solution space $\mathcal{P}_M(\Omega^d)$ as follows:
\begin{equation*}
	L:= \begin{cases}
		(M,N), \ &d=1 \\
		(M,M,N), \ &d=2
	\end{cases}, 
    \quad
    \mathcal{P}_M(\Omega^d):= \mathcal{P}_M(\Omega)^d.
\end{equation*}
In the M\"untz spectral method for solving equation \eqref{TFDE}, we seek a numerical solution $u_L(\boldsymbol{x},t)\in \mathcal{P}_M(\Omega^d)\otimes S^{\alpha,-1}_{N,\lambda}(I)$ such that, for any test function $v\in \mathcal{P}_M(\Omega^d)\otimes S^{\alpha,-1}_{N,\lambda}(I)$,
\begin{equation}\label{CalDouLi}
	\mathcal{A}(u_L,v)=\mathcal{F}(v),
\end{equation}
where the bilinear form $\mathcal{A}(u_L,v)$ and the functional $\mathcal{F}(v)$ are defined by
\begin{align*}
	&\mathcal{A}\left(u_L,v\right):=\left(D^{\mu}_t u_L,v\right)_{L^2(Q)}  +\kappa\left(\nabla u_L,\nabla v\right)_{L^2(Q)}+\rho \left(\nabla u_L,v\right)_{L^2(Q)},\\
	&\mathcal{F}\left(v\right):=\left(f,v\right)_{L^2(Q)}.
\end{align*}

\subsection{Efficient implementation}
An efficient implementation of the M\"untz spectral scheme \eqref{CalDouLi} is presented in this subsection.
We start with the one-dimensional case, where the numerical solution is
\begin{equation}\label{uh}
	u_L(x,t)=\sum^{M-2}_{m=0}\sum^N_{n=0}\hat{u}_{nm}\phi_m(x)J^{\alpha,-1,\lambda}_{n+1}(t),
\end{equation}
with $\left\{\hat{u}_{nm}\right\}$ being the coefficients to be determined. By substituting \eqref{uh} to \eqref{CalDouLi} and taking $v=\phi_p(x) J^{\alpha,-1,\lambda}_{q+1}(t)$, we obtain
\begin{align}\label{matrix_re}
	\begin{split}
		\sum^{M-2}_{m=0}\sum^N_{n=0}\hat{u}_{nm}\Big\{\left(\phi_m,\phi_p\right)\left({}^CD^{\mu}_tJ^{\alpha,-1,\lambda}_{n+1},J^{\alpha,-1,\lambda}_{q+1}\right)+\kappa
		\left(\phi'_m,\phi'_p\right)\left(J^{\alpha,-1,\lambda}_{n+1},J^{\alpha,-1,\lambda}_{q+1}\right)\\+\rho\left(\partial_x\phi_m,\phi_p\right)\left(J^{\alpha,-1,\lambda}_{n+1},J^{\alpha,-1,\lambda}_{q+1}\right)\Big\}
		=\left(f,\phi_pJ^{\alpha,-1,\lambda}_{q+1}\right).
	\end{split}
\end{align}
It is noted that the inner products $(\cdot,\cdot)_{L^2(\Omega)}$, $(\cdot,\cdot)_{L^2(I)}$ and $(\cdot,\cdot)_{L^2(I,L^2(\Omega))}$ are simply denoted as $(\cdot,\cdot)$ here.
By defining the following vectors and matrices
\begin{align*}
	&\mbox{U}=\left(\hat{u}_{nm}\right)_{0\leq n\leq N,\ 0\leq m\leq M-2},\\
	&\mbox{F}=(f_{nm})_{0\leq n\leq N,\ 0\leq m\leq M-2},\\
	&\mbox{S}^x=(s^x_{jk})_{0\leq j,k\leq M-2},~\mbox{M}^x=(m^x_{jk})_{0\leq j,k\leq M-2},~\mbox{C}^x=(c^x_{jk})_{0\leq j,k\leq M-2},\\
	&\mbox{S}^t=\left(s^t_{qn}\right)_{0\leq q,n\leq N},\ \mbox{M}^t=\left(m^t_{qn}\right)_{0\leq q,n\leq N},
\end{align*}
where
\begin{align*}
	&f_{nm}=\left(f,\phi_m(x)J^{\alpha,-1,\lambda}_{n+1}\left(t\right)\right),\\
	&s^t_{qn}=\left(^CD^{\mu}_tJ^{\alpha,-1,\lambda}_{n+1}(t),~J^{\alpha,-1,\lambda}_{q+1}(t)\right),\\
	&m^t_{qn}=\left(J^{\alpha,-1,\lambda}_{n+1}(t),J^{\alpha,-1,\lambda}_{q+1}(t)\right),
\end{align*}
numerical scheme \eqref{matrix_re} can be rewritten in a compact matrix form
\begin{equation}\label{Mform_1d}
	\mbox{S}^t\mbox{U}\mbox{M}^x+\kappa\mbox{M}^t\mbox{U}\mbox{S}^x+\rho\mbox{M}^t\mbox{U}\mbox{C}^x = \mbox{F}.
\end{equation}
The above matrix system is not large in general and can be solved efficiently by using a direct method.

However, the space-time spectral method for a multi-dimensional problem is computationally expensive, suffering from the curse of dimensionality in the sense that the size of equation system grows quickly with the dimension of space.
Therefore, for two and higher dimensional cases, it is necessary to develop an efficient solver for the space-time M\"untz spectral method.

For a two-dimensional sub-diffusion problem described by \eqref{TFDE} with $\rho=0$, the M\"untz spectral method can be efficiently implemented by using Fourier-like basis functions in space \cite{Shen2007Fourierization}.
The mass matrix in space can be diagonalized through eigendecomposition. For the mass matrix $\mbox{M}^x$ in the one-dimensional case, the eigendecomposition is
\begin{equation*}
	\mbox{M}^x\mbox{E}=\mbox{E}{\Lambda},~~~\mbox{E}^{\top}\mbox{E}=\mbox{I}_d,
\end{equation*}
where ${\Lambda}=\mbox{diag}(\bar\lambda_j)$ is the diagonal matrix whose diagonal elements are the eigenvalues, $\mbox{E}=(e_{jk})_{j,k=0,1,\dots,M-2}$ is the normalized eigenvector matrix, and $\mbox{I}_d$ is the identity matrix.
Then we define functions
\begin{equation*}
	\sigma_k(x):=\sum^{M-2}_{j=0}e_{jk}\phi_j(x), k=0,1,\dots,M-2,
\end{equation*}
where $\{\phi_j\}$ are the basis functions defined in \eqref{basis}.
The functions $\{\sigma_k(x)\}$ form a new basis of $\mathcal{P}_M$ satisfying
\begin{equation*}
	\int^1_{-1} \sigma_k\sigma_l dx=\bar\lambda_k\delta_{kl},~~~\int^1_{-1}\sigma^{\prime}_k\sigma^{\prime}_ldx=\delta_{kl}.
\end{equation*}
By using this basis in space, the numerical solution of the two-dimensional sub-diffusion problem can be expressed as
\begin{equation}\label{sol_exp}
	u_L(x,y,t)
	=\sum^N_{n=0} \sum^{M-2}_{k,l=0}\hat{u}_{k,l,n}\sigma_k(x)\sigma_l(y)J^{\alpha,-1,\lambda}_{n+1}(t),
\end{equation}
where $\left\{\hat{u}_{k,l,n}\right\}$ are the coefficients to be determined. The mass and stiffness matrices in space are diagonalized as
\begin{equation*}
	\mbox{M}^{x,y}=\mbox{diag}(\bar\lambda_k)\otimes\mbox{diag}(\bar\lambda_k),~~\mbox{S}^{x,y}=\mbox{I}_d\otimes \mbox{diag}(\bar\lambda_k)+\mbox{diag}(\bar\lambda_k)\otimes \mbox{I}_d,
\end{equation*}
taking advantage of the Fourier-like property of the new basis $\{\sigma_k(x)\}$.

It is known that the eigendecomposition of a non-symmetric matrix is numerically unstable \cite{Shen2019efficient}. Consequently, instead of the aforementioned eigendecomposition, a QZ decomposition is applied to the matrices in time due to the non-symmetric nature of the stiffness matrix $\mbox{S}^t$.
It is noted that the QZ decomposition has essentially the same computational complexity as the eigendecomposition \cite{Shen2019efficient,chen2020spectrally}.
The QZ decomposition is
\begin{equation}\label{qz}
	\mbox{Q}({\mbox{S}^t})^\top \mbox{Z}={\mbox{A}},~~\mbox{Q}({\mbox{M}^t})^\top \mbox{Z}={\mbox{B}},
\end{equation}
where ${\mbox{A}}=({a}_{ij})$ and ${\mbox{B}}=({b}_{ij})$ are upper triangular matrices, $\mbox{Q}$ and $\mbox{Z}$ are orthogonal matrices satisfying
\begin{equation*}
	\mbox{Q}\mbox{Q}^\top=\mbox{I}_d, \quad \mbox{Z}\mbox{Z}^\top=\mbox{I}_d.
\end{equation*}

The matrix form of the M\"untz spectral scheme for a two-dimensional sub-diffusion problem is 
\begin{equation}\label{2D_system_matrix}
	\mbox{S}^t\mbox{U}(\mbox{M}^{x,y})^\top+\kappa\mbox{M}^t \mbox{U} (\mbox{S}^{x,y})^\top=\mbox{F},
\end{equation}
where
\begin{equation*}
	\begin{gathered}
		\mbox{U}=
		\begin{bmatrix}
			\hat{u}_{0,0,0} & \ldots & \hat{u}_{0,M-2,0}&\ldots&\hat{u}_{M-2,0,0}&\ldots&\hat{u}_{M-2,M-2,0}\\
			\vdots& &\vdots& &\vdots& &\vdots\\
			\hat{u}_{0,0,N}&\ldots&\hat{u}_{0,M-2,N}&\ldots&\hat{u}_{M-2,0,N}&\ldots&\hat{u}_{M-2,M-2,N}
		\end{bmatrix},
	\end{gathered}
\end{equation*}
is the matrix of unknown coefficients.
By substituting \eqref{qz} into \eqref{2D_system_matrix}, and then multiplying \eqref{2D_system_matrix} by $\mbox{Z}^\top$, we obtain the following matrix system
\begin{equation}\label{2D_system_efficient}
	{\mbox{A}}^\top \mbox{V}(\mbox{M}^{x,y})^\top+\kappa{\mbox{B}}^\top \mbox{V} (\mbox{S}^{x,y})^\top=\mbox{Z}^\top{\mbox{F}},
\end{equation}
with $\mbox{U}=\mbox{Q}^\top \mbox{V}$.
The matrix system \eqref{2D_system_efficient} can be solved efficiently, due to the simple structures of its matrices.
It is worth noting that the computational cost of solving the above system is equivalent to that of solving $N+1$ two-dimensional elliptic problems, see \cite{Shen2019efficient,Zeng2024An} for further technical details.

\begin{remark}
	Computation of the entries of matrices $S^t$ and $M^t$ is crucial for implementation of the M\"untz spectral method.
	An efficient evaluation method proposed in \cite{Hou2018Muntz} is
	\begin{align*}
		&s^t_{qn}=\frac{(n+\alpha+1)(q+\alpha+1)}{\lambda\Gamma(1-\mu)(n+1)}\sum^N_{i=0}\sum^{\tilde{N}}_{j=0}
		(\frac{1-\hat{\zeta}^{1/{\lambda}}_j}{1-\hat{\zeta}_j})^{-\mu}J^{\alpha+1,0}_n(2\zeta_i\hat{\zeta}_j-1)\hat{\omega}_j
		J^{\alpha,1}_q(2\zeta_i-1)\omega_i, \\
		& m^t_{qn}=\frac{(n+\alpha+1)(q+\alpha+1)}{\lambda(n+1)(q+1)}\sum^N_{i=0}\rho_iJ^{\alpha,1}_n(2\xi_i-1)J^{\alpha,1}_q(2\xi_i-1),
	\end{align*}
	where the Gauss quadrature point sets $\{\zeta_i\}^N_{i=0}$, $\{\hat{\zeta}_j\}^{\hat{N}}_{j=0}$ and $\{\xi_i\}^N_{i=0}$ are the zeros of shifted Jacobi polynomials $J^{0,\frac{1-\mu}{\lambda}+1}_{N+1}(2t-1),$ $J^{-\mu,0}_{\hat{N}+1}(2t-1)$ and $J^{0,1+1/{\lambda}}_{N+1}(2t-1)$, respectively, $\{\omega_m\}^N_{m=0},$ $\{\hat{\omega}_n\}^{\hat{N}}_{n=0}$ and $\{\rho_i\}^N_{i=0}$ are the associated weights.
\end{remark}

\subsection{Empirical parameter selection}\label{issue}

In the M\"untz spectral method, $\lambda$ is a free parameter that lies in the range of $(0,1]$. Different values of $\lambda$ results in different M\"untz Jacobi polynomial bases, and hence different numerical solutions. Therefore, the choice of $\lambda$ is crucial to the accuracy of the M\"untz spectral method. 
However, there is no theoretical guidance or reliable rule for determining the optimal values of $\lambda$, making its selection a significant challenge in practical implementation.

The one-dimensional sub-diffusion equation, described by \eqref{TFDE} with $\rho=0$ and $\kappa=1$, can be used as an example to illustrate the challenges in parameter selection. Under the assumption that the forcing function $f$ exhibits certain regularity, a solution representation \cite{sakamoto2011initial,jin2013error} can be obtained in terms of Mittag-Leffler function for the sub-diffusion equation as follows
\begin{align}\label{Sol_ML}
	\begin{split}
		u(x, t) & =\sum_{i=1}^{\infty}\left[\int_0^t\left(f(\cdot, \tau), \psi_i(\cdot)\right)(t-\tau)^{\mu-1} E_{\mu, \mu}\left(-\bar \lambda_i(t-\tau)^\mu\right) d \tau\right] \psi_i(x) \\
		& =t^\mu \sum_{i=1}^{\infty}\left[\int_0^1\left(f(\cdot, \tau t), \psi_i(\cdot)\right)(1-\tau)^{\mu-1} E_{\mu, \mu}\left(-\bar \lambda_i t^\mu(1-\tau)^\mu\right) d \tau\right] \psi_i(x),
	\end{split}
\end{align}
where $-\partial^2_{xx}\psi_i(x) =\bar \lambda_i\psi_i(x)$, $\psi_i(\pm 1)=0$, with $\{\bar{\lambda}_i\}^{\infty}_{i=0}$ and $\{\psi_i(x)\}^{\infty}_{i=0}$ being the eigenvalues and eigenfunctions of the operator $-\partial^2_{xx}$ on the domain $\Omega$ subject to the homogeneous Dirichlet condition, respectively, and the Mittag-Leffler function $E_{\mu,\nu}(z)$ is defined by
\begin{equation*}
	E_{\mu,\nu}(z)=\sum_{n=0}^\infty \frac{z^n}{\Gamma(\mu n+\nu)},\quad \mu>0,\ \nu\in\mathbb{R},\ z\in \mathbb{C}.
\end{equation*}
For more details on the Mittag-Leffler function, see \cite{Podlubny1998Fractional}. It is observed from \eqref{Sol_ML} that, even given a smooth forcing function $f$, the solution $u$ may exhibit
singularity with the leading order $t^\mu$ at the starting point $t = 0$.

The choice of $\lambda$ is of great importance for the
accuracy of the M\"untz spectral method in solving the sub-diffusion problem, yet there is no systematic strategy for selecting optimal parameter values.
While empirical values of $\lambda$ are suggested in the literature for specific cases, these recommendations lack theoretical justification. In practice, researchers often test multiple values of $\lambda$ based on prior experience, choosing the most effective option. 
We can discuss the parameter selection issue case by case.

Case I: if the solution $u$ is known to be smooth, there is no doubt that the optimal value is $\lambda=1$. In this case, the M\"untz spectral method reduces to the classical spectral method based on the usual Jacobi polynomials.

Case II: if the forcing function $f(x,t)$ is smooth and $\mu=p/q$ is a rational number, where $p$ and $q$ are positive integers, a common practice is to set $\lambda=1/q$. Doing so, the mapped solution $u(x,t^{1/\lambda})$ attains sufficient smoothness, enabling the spectral method to achieve spectral accuracy. However, there is a defect of this choice in the specific implementation. The matrix system of the M\"untz spectral method depends on both $\lambda$ and $\mu$. In practice, we observe that when $\lambda$ is small, the condition number of the matrix system becomes very large, leading to numerical instability and inaccurate solutions or even simulation failures. 
One cause of this issue lies in the definition of the differential operators in fractional Jacobi-weighted Sobolev spaces:
\begin{equation*}
	D^0_{\lambda}:=I_d,~~D^1_{\lambda}:=\frac{d}{dx^{\lambda}}=\frac{d}{\lambda x^{\lambda-1}dx},~~D^k_{\lambda}:=\overbrace{D^1_{\lambda}D^1_{\lambda}\cdots D^1_{\lambda}}^k.
\end{equation*}
It is observed that each time the derivative $D^1_{\lambda}$ is taken, the denominator will be multiplied by $\lambda$. When $q$ is large, $\lambda$ is small, introducing singularity to the system.
This analysis implies that the choice $\lambda=1/q$ is invalid when $\mu$ is close to zero.

Case III: if the forcing function $f(x,t)$ is smooth and $\mu$ is an irrational number, there is no suitable value of $\lambda$ to make the mapped solution $u(x,t^{1/\lambda})$ sufficiently smooth. In this case, we can set $\lambda=1/q$ with a reasonably large $q$ such that $u(x,t^{1/\lambda})$ is smooth enough. The integer $q$ should be large to improve the smoothness of the mapped solution, while not overly large to avoid numerical instability, as discussed in Case II. Therefore, it is difficult to select a proper value of $\lambda$ in this case.

In summary, the principle of selecting $\lambda$ is to ensure that $u(x,t^{1/\lambda})$ is sufficiently smooth, which seems to be straightforward. However, in practical implementation, selecting a parameter value that yields an accurate numerical solution proves challenging due to the absence of rigorous theoretical guidance.
As a result, $\lambda$ is often chosen empirically by testing multiple values and selecting the one that yields the highest accuracy \cite{Hou2018Muntz,Hou2017A}. This time-consuming parameter tuning process has become a major bottleneck in applying the M\"untz spectral method.

\section{Machine learning optimization of M\"untz spectral method}
\label{sec3}

This section presents an optimization approach for the M\"untz spectral method based on machine learning. It is pointed out in Section \ref{sec2} that, there is a free parameter $\lambda$ that affects the performance of the M\"untz spectral method. Improper choice of $\lambda$ results in inaccurate solution or numerical instability. 
In practical implementations, values of $\lambda$ are determined empirically, as no established theoretical framework is available to guide the selection process. A major challenge in applying the M\"untz spectral method is the time-consuming process of parameter tuning.

A machine learning prediction model of $\lambda$ is developed in this section to address the parameter selection issue.  An artificial neural network (ANN) is used to predict the optimal value of $\lambda$, given the time-fractional order $\mu$. The optimal value of the parameter $\lambda$ can be pre-determined by the trained ANN, which eliminates any additional computational overhead for the M\"untz spectral solver. 
The structure and training of the ANN will be presented in detail in the remainder of this section.

\subsection{Artificial neural network}\label{subsec:feedforward_network}
For the M\"untz spectral method, the optimal value of the parameter $\lambda$, which is used to compute the fractional polynomial basis functions in the time direction, depends on the solution regularity that is determined by the governing equation and the initial condition. Given that we focus on a homogeneous initial condition, the time-fractional order $\mu$ and the forcing function $f$ are the key factors influencing the optimal value of $\lambda$. In this paper, we train a machine learning model to approximate the mapping from the time-fractional order $\mu$ to the optimal value of $\lambda$, with the aim of applying this model to general time-fractional PDEs with varying forcing functions.

\begin{figure}[htbp!]
	\centering
	\includegraphics[width=0.7\linewidth]{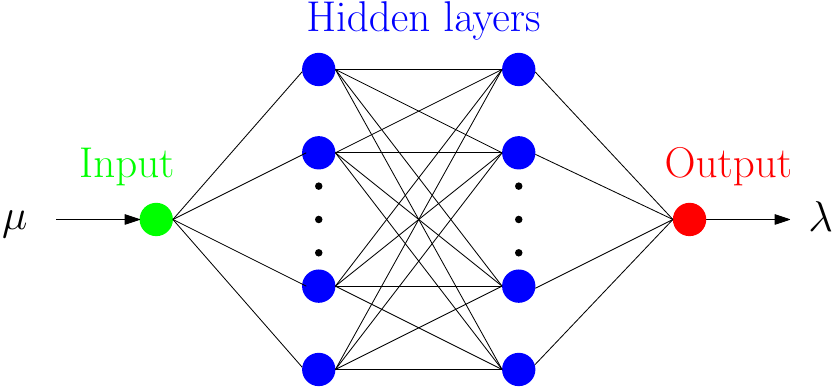}
	\caption{A sample feedforward neural network structure with two hidden layers.}
	\label{fig:network_structure}
\end{figure}

A feedforward neural network (FNN) is employed to predict the optimal value of the parameter $\lambda$, given the time-fractional order $\mu$. The FNN consists of an input layer, $L$ hidden layers and an output layer, as shown in Fig. \ref{fig:network_structure}. The $l$-th layer has $n_l$ artificial neurons, where $l=0,1, \cdots,L, L+1$ denote the input layer, first hidden layer, $\cdots$, $L$-th hidden layer and the output layer, respectively. In this work, the input dimension is $n_0=1$ and the output dimension is $n_{L+1}=1$. The hidden layers are set to have the same number of neurons, i.e., $n_1=n_2=\cdots=n_L=n_H$. 

The network output $\lambda$ is a nonlinear function of the input $\mu$. The evaluation of the function $\lambda=f_{\mathrm{NN}} \left(\mu\right)$, which is the forward propagation of the network,  can be expressed as
\begin{equation}\label{network-propagation}
	\begin{aligned}
		&\mbox{y}^{\left(0\right)}= \mu, \\
		& \mbox{y}^{\left(l\right)}= f^{\left(l\right)}_{\mathrm{act}} \left(\mbox{W}^{\left(l\right)} \mbox{y}^{\left(l-1\right)} + \mbox{b}^{\left(l\right)}\right), \quad l=1,\cdots, L+1, \\
		& \lambda= \mbox{y}^{\left(L+1\right)},
	\end{aligned}
\end{equation}
where $\mbox{W}^{\left(l\right)} \in \mathbb{R}^{n_l \times n_{l-1}}$, $\mbox{b}^{\left(l\right)} \in \mathbb{R}^{n_l }$ and $\mbox{y}^{\left(l\right)} \in \mathbb{R}^{n_l}$ are the weight, bias and output of the $l$-th layer, respectively. $f^{\left(l\right)}_{\mathrm{act}}$ is the element-wise activation function of the $l$-th layer. In this work, the activation function of each hidden layer is the Leaky ReLU function 
\begin{equation*}
	f^{\left(l\right)}_{\mathrm{act}} \left(x\right)=\max\left(0,x\right) + {negative\_slope}*\min\left(0,x\right), \quad l=1, \cdots,L,
\end{equation*}
where the parameter $negative\_slope$ is set to $0.01$, following the default configuration in PyTorch \cite{paszke2019pytorch}.
The Sigmoid function
\begin{equation*}
	f^{\left(L+1\right)}_{\mathrm{act}} \left(x\right)= \frac{1}{1+e^{-x}},
\end{equation*}
is used as the activation function of the output layer to ensure that the predicted $\lambda$ stays in the range $\left(0,1\right)$.

The weights $\{ \mbox{W}^{\left(l\right)} \}^{L+1}_{l=1}$ and biases $\{ \mbox{b}^{\left(l\right)} \}^{L+1}_{l=1}$ are adjustable parameters that are optimized during the training process described in Section \ref{subsec:network_training}. The number of hidden layers $L$ and the number of neurons in each hidden layer $n_H$ are hyper-parameters that can be tuned to obtain an accurate and efficient network configuration.

\subsection{Supervised learning}\label{subsec:optimization_strategy}

In this subsection, we present an optimization strategy for training a FNN to predict the optimal value of $\lambda$. As discussed in Section \ref{sec2}, varying $\lambda$ results in different M\"untz Jacobi polynomial bases in the time domain, which in turn leads to different numerical solutions. The optimal $\lambda$ is the one that yields the most accurate space-time numerical solution. 

A suitable value of the parameter $\lambda$ is selected to construct an appropriate M\"untz Jacobi polynomial basis, effectively capturing the singularity of the solution in time. For a given time-fractional order $\mu$, the optimal value of $\lambda$ is assumed to perform consistently across different spatial dimensions. This allows for efficient training of a FNN on a one-dimensional spatial domain, with the resulting model being directly applicable to higher-dimensional problems.

Based on the above analysis, we propose to train a FNN as the optimal parameter prediction model for the M\"untz spectral method, by minimizing the solution errors of the following one-dimensional time-fractional convection-diffusion equation
\begin{equation}\label{model_equation}
	D_t^\mu u(x,t)= \Delta u(x,t)-\nabla u(x,t)+f(x,t),\ x \in \Omega,~t\in I.
\end{equation}
Given homogeneous initial and boundary conditions, the exact solution of \eqref{model_equation} is determined by the forcing function $f\left(x,t\right)$. With general forcing functions, the exact solutions are not available and the solution errors can not be computed directly. To address this issue, a manufactured solution
\begin{equation}\label{manufactured_solution}
	u_e\left(x,t\right)= t^\nu \sin(\pi x), \quad 0<\nu<1,
\end{equation}
is employed as the exact solution of \eqref{model_equation}, with the corresponding forcing function derived analytically by balancing the equation. The parameter $\nu$ is used to control the singularity of the exact solution.

\begin{figure}[tb!]
	\centering
	\includegraphics[width=\linewidth]{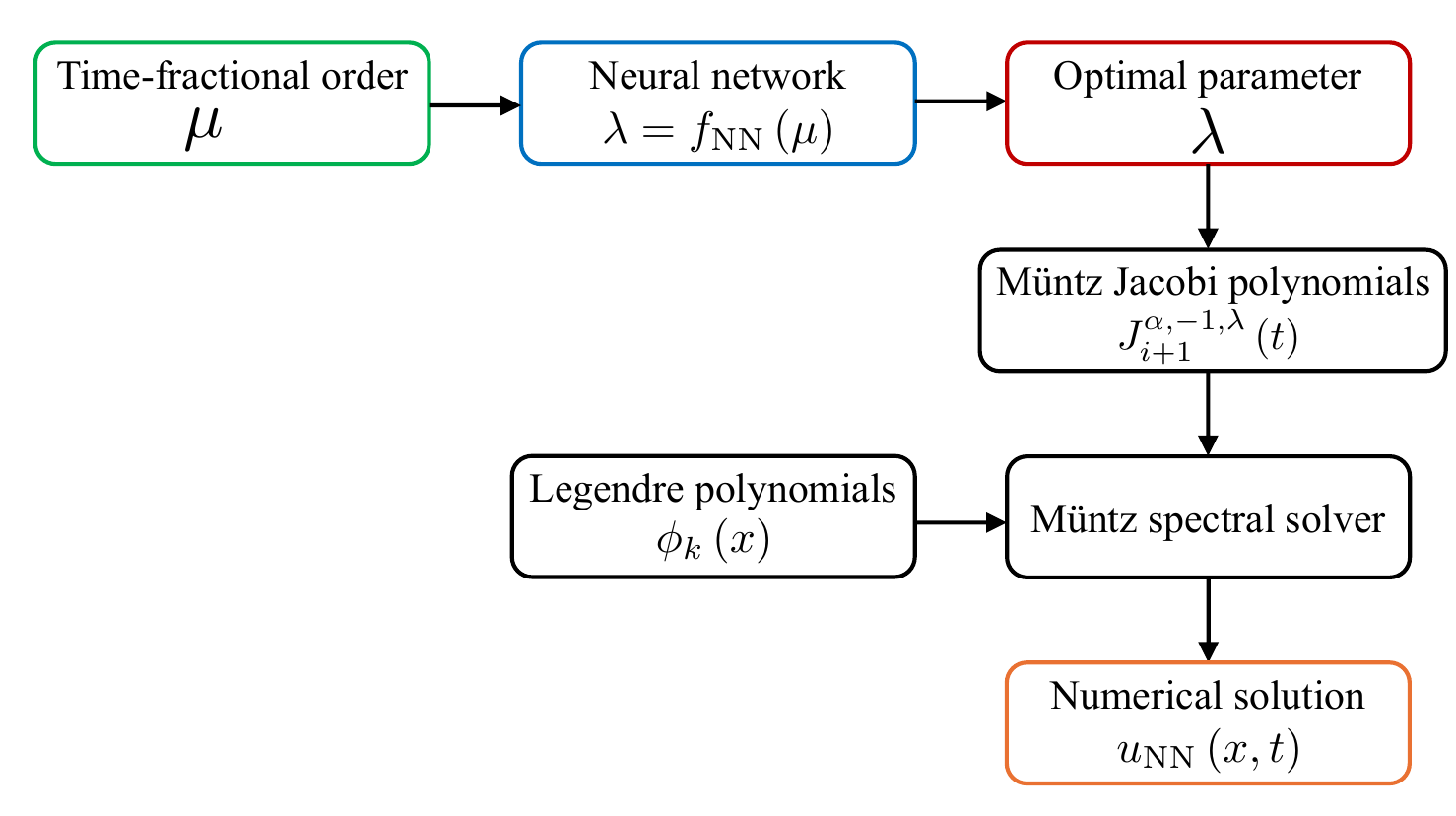}
	\caption{Workflow to compute the numerical solution.}
	\label{fig:procedures}
\end{figure}

A numerical solution $u_{\mathrm{NN}}$ for the equation \eqref{model_equation} can be obtained by using the M\"untz spectral method, with a FNN as the parameter prediction model, following the workflow illustrated in Fig. \ref{fig:procedures}. 
The error of the numerical solution is computed as
\begin{equation*}
	E \left(u_e, u_{\mathrm{NN}}\right)= \left\| \left( u_e \left(x_i, t_j \right)  - u_{\mathrm{NN}} \left(x_i, t_j \right) \right)_{0\leq i \leq M-2, \ 0 \leq j \leq N} \right\|_{\mathrm{F}},
\end{equation*}
where $x_i$ and $t_j$ are the coordinates of the Gauss quadrature points in space and time, respectively.

The FNN is trained using a supervised learning paradigm, where the input feature is the parameter instance $\left(\mu, \nu\right)$, and the label is the exact space-time solution $u_e\left(x,t;\nu\right)$. The optimal network configuration is obtained by minimizing the discrepancy between the numerical solution $u_{\mathrm{NN}}$ and the exact solution $u_e$, specifically the solution error $E \left(u_e, u_{\mathrm{NN}}\right)$, throughout the training process. 

It is noted that, in the M\"untz spectral method, the parameter $\lambda$ only influences the M\"untz Jacobi polynomial basis in the time direction. Therefore, a sufficient number of polynomial basis functions are required in the space direction to ensure that the solution error predominantly arises from the time component, facilitating effective optimization of the parameter prediction model. In this work, the numbers of basis functions in time and space are $N+1=11$ and $M-1=19$, respectively.

\subsection{Dataset generation}\label{subsec:dataset_generation}

This subsection presents the generation of datasets used to train the FNN. A training dataset is denoted as $\mathcal{D}_{tr}= \left\{\left(\mu_i, \nu_i\right), u_e(x,t;\nu_i)\right\}^{N_{tr}}_{i=1}$, where $N_{tr}$ is the dataset size. The training dataset is generated by first sampling a set of parameter instances on the two-dimensional domain $\mathcal{P}_{\mu,\nu} = \left(0,1\right) \times \left(0,1\right)$, and then computing the corresponding exact solutions according to \eqref{manufactured_solution}. 

In this work, a training dataset of size $N_{tr}= 900$ is constructed using tensorial (or full factorial) sampling \cite{quarteroni2015reduced}, incorporating $30$ randomly sampled points along the $\mu$-axis and $30$ randomly sampled points along the $\nu$-axis. 

In supervised learning, a validation dataset is often used to prevent over-fitting. In this work,  a validation dataset $\mathcal{D}_{va}= \left\{\left(\mu_i, \nu_i\right), u_e(x,t;\nu_i)\right\}^{N_{va}}_{i=1}$ of size $N_{va}= 900$ is generated using tensorial sampling, incorporating $30$ uniformly sampled points along the $\mu$-axis and $30$ uniformly sampled points along the $\nu$-axis. 

\subsection{Training process}\label{subsec:network_training}

%

The training of the network is implemented in PyTorch \cite{paszke2019pytorch}. The optimal weights and biases of the network are obtained using the Adam stochastic optimizer \cite{kingma2014adam}, which uses mini-batches of size $N_b= 100$ of the training data set $\mathcal{D}_{tr}$ to take a single step by minimizing the loss function
\begin{equation}\label{eqn:Loss_function}
	\mathcal{L}= \dfrac{1}{N_b} \sum^{N_b}_{i=1}
	\dfrac{ E \left(u_e \left(x,t;\nu_i\right), u_{\mathrm{NN}} \left(x,t;\mu_i, \nu_i\right) \right) } { E \left(u_e \left(x,t;\nu_i\right), u_{\mathrm{REF}} \left(x,t;\mu_i, \nu_i\right) \right) }.
\end{equation} 
The reference solution $u_{\mathrm{REF}}$ is computed using the M\"untz spectral method with $\lambda=1$, which corresponds to the classical spectral method based on the usual Jacobi polynomial basis functions.
The solution error is normalized by the error of the reference solution, as it is observed in practice that solution errors across different cases can vary significantly in magnitude. This normalization ensures better control over the relative importance of solution errors for different singularities.

The training is performed for $400$ epochs to obtain a converged network.
The convergence speed is controlled by a learning rate $\eta$, which is dynamically adjusted using the Cosine Annealing Warm Restarts (CAWR) scheduler in PyTorch. 
This scheduler employs a strategy that combines cosine annealing with periodic learning rate restarts, which has been shown to facilitate the escape from local minima and potentially enhance convergence to a global optimum. 
The training process is summarized in Algorithm \ref{alg:training-FNN}.

The training of the network is performed $10$ times, following a multiple restarts approach \cite{hsu1995artificial}, to prevent the training results from depending on the random initialization of the weights and biases. The selected trained network is the one that has the smallest loss on the validation dataset $\mathcal{D}_{va}$.

Two important hyper parameters for the FNN used in this work are the number of hidden layers $L$ and the number of neurons in each hidden layer $n_H$, which determine the structure of the network. In practice, we progressively increased the network size and found that the optimal configuration, balancing accuracy and sparsity, was achieved with $L = 2$ hidden layers and $n_H$ = 20 neurons per layer.

\begin{algorithm}[htbp!]
	\caption{Training process of the FNN.}
	\label{alg:training-FNN}
	\begin{algorithmic}[1]
		
		%
		
		\Function{[$\mathbf{W}_{tr}, \mathbf{b}_{tr}$]=FNN\_TRAINING}{$\mathcal{D}_{tr}, N_{epo}, N_{b}$}
		\State $\mathbf{W}, \mathbf{b} \gets \mathrm{INIT}(\mathbf{W},\mathbf{b})$ \Comment {weight and bias initialization}
		\State $ N_{batch} \gets N_{tr}/N_{b}$\Comment{number of mini-batches}
		\For{$epoch \gets 1, N_{epoch}$} \Comment{training epoch loop}
		\State $\eta \gets \mathrm{CAWR}(epoch)$  \Comment{CAWR learning rate scheduler}
		\State $\mathcal{D}_{tr} \gets \mathrm{SHUFFLE}(\mathcal{D}_{tr})$ \Comment{training dataset shuffling}
		\For{$s \gets 1, N_{batch}$} \Comment{mini-batch loop}
		\State $\mathcal{D}_{tr,s} \gets \mathcal{D}_{tr}[(s-1)*N_b +1 :s*N_b,:]$ \Comment{the $s$-th mini-batch}
		\vspace{1mm}
		\State $\mathcal{L} \gets \dfrac{1}{N_b} \sum^{N_b}_{i=1}
		\dfrac{ E \left(u_e \left(x,t;\nu_i\right), u_{\mathrm{NN}} \left(x,t;\mu_i, \nu_i\right) \right) } { E \left(u_e \left(x,t;\nu_i\right), u_{\mathrm{REF}} \left(x,t;\mu_i, \nu_i\right) \right) }$ \Comment{loss}
		\vspace{1mm}
		\State $\Delta \mathbf{W} \gets - \eta \nabla_{\mathbf{W}} \mathcal{L}, \ \Delta \mathbf{b} \gets - \eta \nabla_{\mathbf{b}} \mathcal{L}$ \Comment{Adam optimizer step}
		\State $\mathbf{W} \gets \mathbf{W} + \Delta \mathbf{W}, \ \mathbf{b} \gets \mathbf{b} + \Delta \mathbf{b}$ \Comment{weight and bias update}
		
		\EndFor
		\EndFor
		
		%
		%
		\State $\mathbf{W}_{tr} \gets \mathbf{W}, \mathbf{b}_{tr} \gets \mathbf{b}$ \Comment{optimized weight and bias}

		\EndFunction

	\end{algorithmic}
\end{algorithm}

%

\begin{figure}[htbp!]
	\centering
	\includegraphics[scale=0.6]{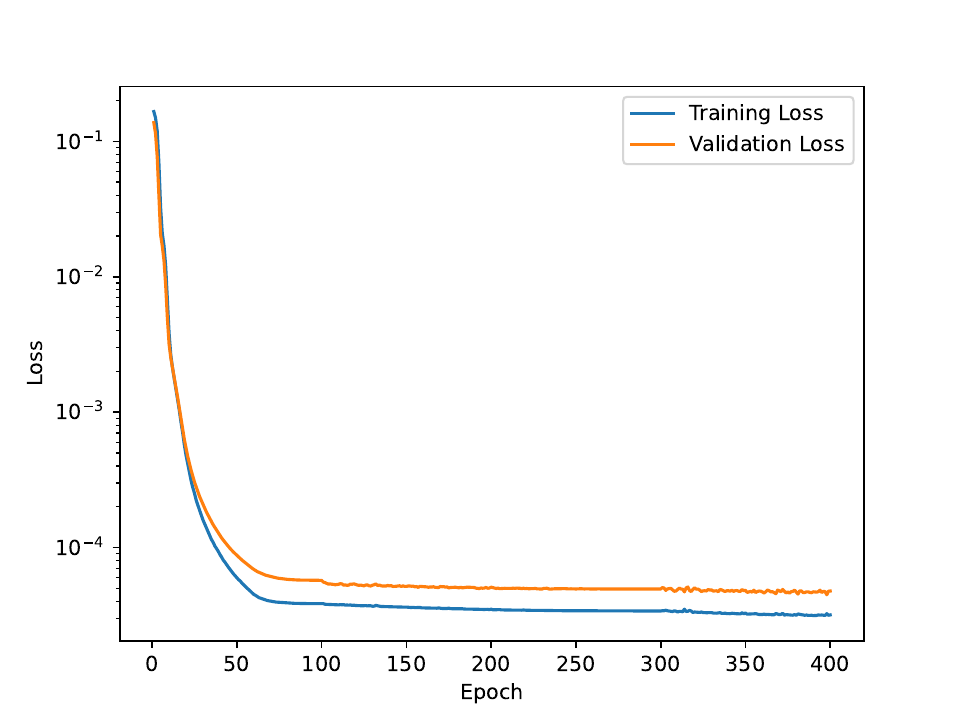}
	\caption{Training and validation loss curves.}
	\label{fig:loss}
\end{figure}

The training and validation loss curves of the final FNN model are shown in Fig. \ref{fig:loss}. It is observed in Fig. \ref{fig:loss} that the converged validation loss is below $10^{-4}$, which means that the M\"untz spectral method using the optimal parameter value predicted by the trained network is at least $4$ orders of magnitude more accurate than the reference method, i.e., the classical spectral method. 

Once the offline training is finished, the weights and biases of the network are saved. The online evaluation of the network is very efficient as it just involves the forward propagation \eqref{network-propagation}. 
Given a new value of $\mu$, the trained network is evaluated to predict the optimal value of $\lambda$ to be used by the M\"untz spectral solver, at a negligible computational cost. Numerical results in Section \ref{sec4} demonstrate that the FNN parameter predictions lead to a substantial improvement in accuracy.

\subsection{Neural network parameter prediction}\label{subsec:prediction_results}

\begin{figure}[htbp!]
	\centering
	\begin{subfigure}{0.48\textwidth}
		\includegraphics[width=\linewidth]{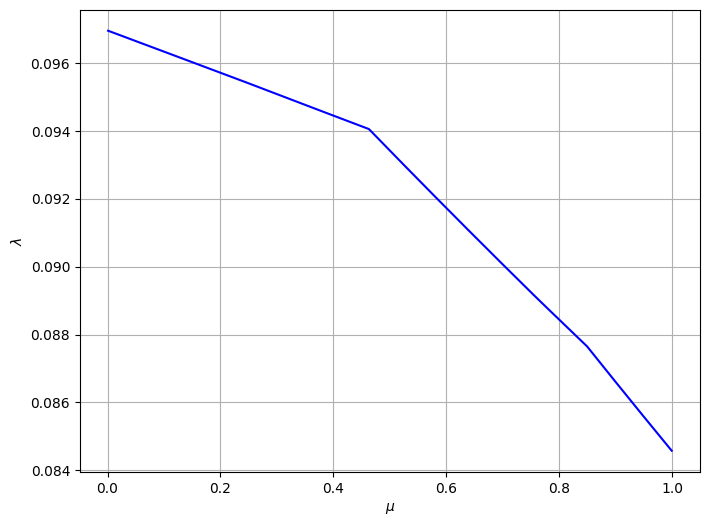}
		\caption{ANN}
		\label{lambda-mu}
	\end{subfigure}
	\begin{subfigure}{0.48\textwidth}
		\includegraphics[width=\linewidth]{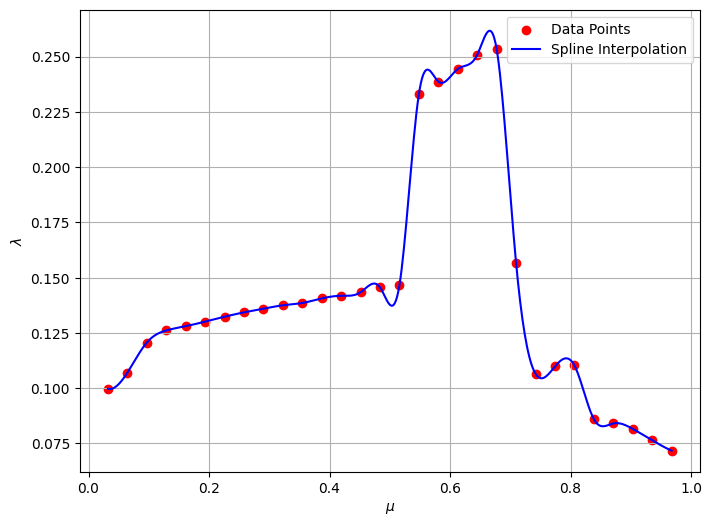}
		\caption{Cubic spline}
		\label{lambda-mu-spline}
	\end{subfigure}
	\caption{Predicted $\lambda$-$\mu$ curves by neural network and cubic spline.}
	\label{prediced-curves}
\end{figure}

The $\lambda$-$\mu$ curve predicted by the trained FNN is shown in Fig. \ref{lambda-mu}. As shown, the FNN predicted $\lambda$ remains within the range $(0.084,0.098)$ and decreases as $\mu$ increases. This trained FNN can be used to predict values of $\lambda$ for the M\"untz spectral method in solving time-fractional convection-diffusion problems, achieving highly accurate numerical solutions. 
It is highlighted that, the FNN trained by minimizing the solution errors of a one-dimensional time-fractional convection-diffusion equation, can be directly applied to multi-dimensional time-fractional PDEs, thus circumventing the need for complex and inefficient model retraining.

\subsection{Comparison with spline interpolation}\label{subsec:Interpolation}


In the previous subsections, we introduced a machine learning-based optimization framework for the M\"untz spectral method. Within this framework, a machine learning model is used to predict the optimal parameter values, trained by minimizing the solution errors of a one-dimensional time-fractional PDE. In principle, a variety of trainable function approximators can serve as the parameter prediction models within this optimization framework. Thus, the machine learning model is not limited to ANNs. For comparison with the ANN, we also train a cubic spline, a classical interpolation model.


The spline interpolation model is constructed through the following three steps:
First, we uniformly sample $30$ values for the time-fractional order $\mu$ within the range $(0, 1)$. Similarly, we sample $30$ values for the parameter $\nu$, which governs the singularity of the manufactured exact solutions, within the range $(0, 1)$.
Second, for each sampled time-fractional order $\mu_i$, we associate $30$ corresponding time-fractional PDEs, each linked to one of the $30$ sampled values of $\nu$. 
Using the Adam optimizer, we determine the optimal value of $\lambda_i$ that minimizes a loss function, defined as the mean normalized absolute solution error across the $30$ PDEs.
Third, after obtaining the dataset $\left\{ (\mu_i, \lambda_i) \right\}_{i=1}^{30}$, we fit a spline to approximate the functional relationship between $\mu$ and the optimal $\lambda$.
In the Adam optimization, the optimal $\lambda$ is sought within the range $[0.001, 0.999]$, with an initial guess $\lambda=0.5$.  The optimization is performed for $400$ iterations. A learning rate scheduler is used to control the convergence speed and an early stopping mechanism is empolyed to enhance efficiency.
The optimization process takes $13,590$ seconds, while training the ANN once requires only $6,315$ seconds. Notably, both processes involve the same number of data points, use the same optimizer, and run for the same number of iterations. These findings demonstrate that the offline training of the ANN is more efficient.


The $\lambda$-$\mu$ curve predicted by the spline interpolation is shown in Fig. \ref{lambda-mu-spline}, which presents the $\left\{ (\mu_i, \lambda_i) \right\}$ data points alongside the cubic spline interpolation model. 
Compared with the curve in Fig. \ref{lambda-mu} predicted by ANN, the spline interpolation spans a broader range and forms a more complex curve.
Notably, for $0.5161 \leq \mu \leq 0.7097$, the spline-predicted $\lambda$ is substantially larger than the corresponding ANN prediction.
Therefore, we construct a numerical example with a time-fractional order $\mu = 0.65 \in (0.5161,0.7097)$  to compare the accuracy of the two parameter prediction models.

\begin{figure}[htbp!]
	\centering
	\begin{subfigure}{0.48\textwidth}
		\includegraphics[width=\linewidth]{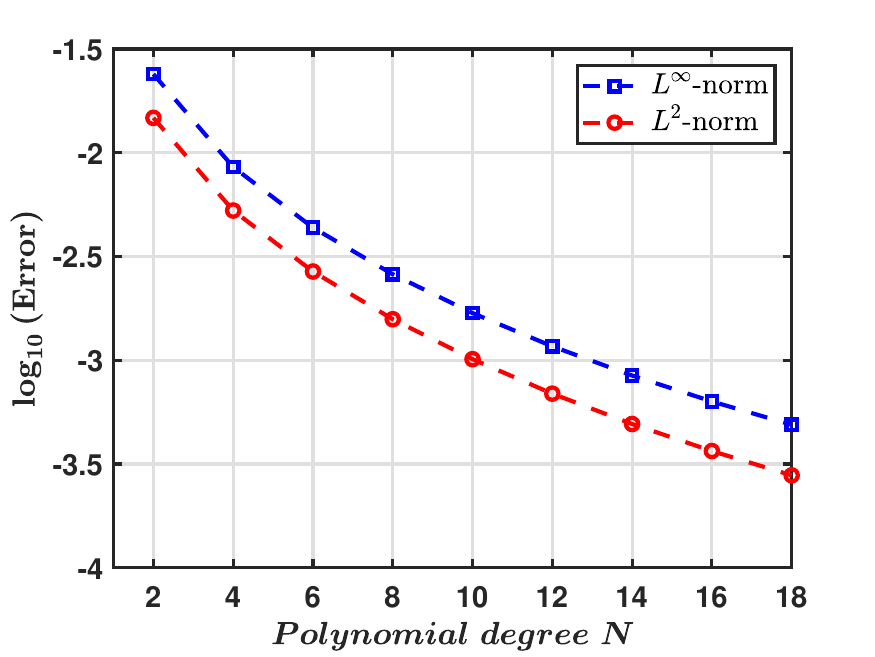}
		\caption{$\mu=0.65$, $\lambda=1$}
		\label{Plot2_1}
	\end{subfigure}
	\begin{subfigure}{0.48\textwidth}
		\includegraphics[width=\linewidth]{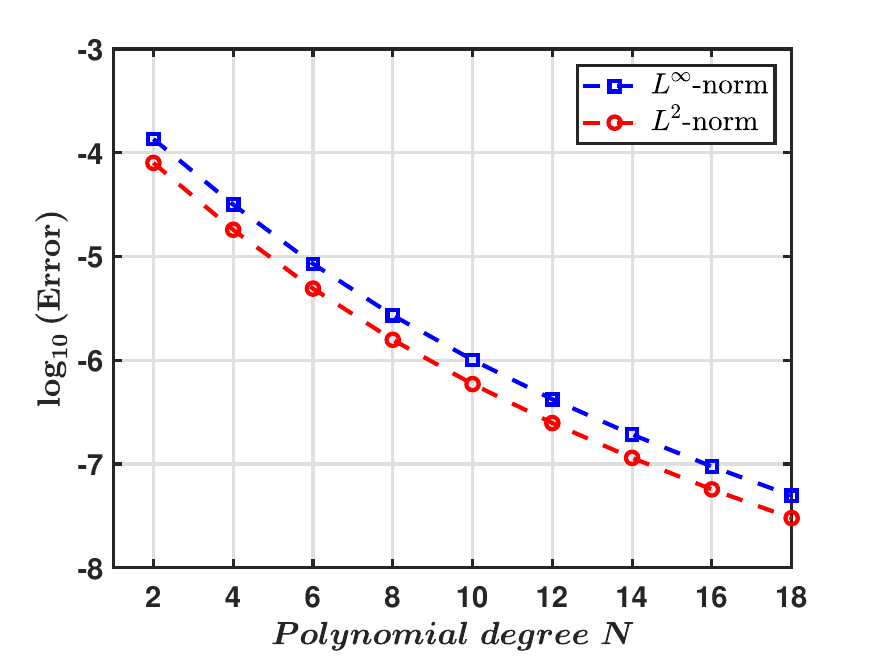}
		\caption{$\mu=0.65$, $\lambda=0.2540$}
		\label{Plot2_2}
	\end{subfigure}
	\begin{subfigure}{0.48\textwidth}
		\includegraphics[width=\linewidth]{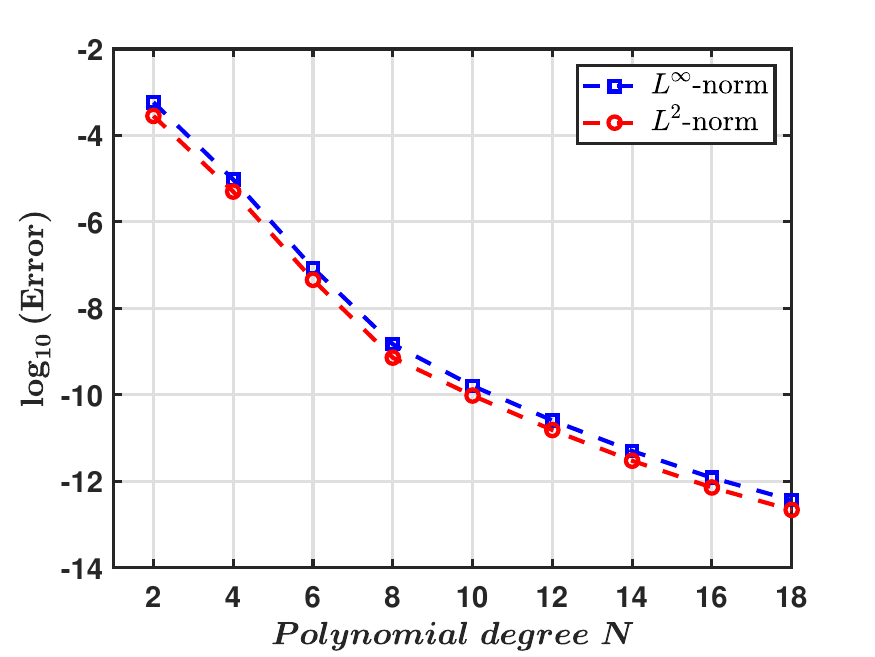}
		\caption{$\mu=0.65$, $\lambda=0.0909$}
		\label{Plot2_3}
	\end{subfigure}
	\caption{(Errors in time) $L^{\infty}$ and $L^2$ errors in semi-log scale versus $N$.}
	\label{Plot2}
\end{figure}

We consider the time-fractional convection-diffusion equation \eqref{TFDE} with a manufactured exact solution $u(x,t)=\sin(\pi x)t^{1-\mu}$. For the time-fractional order $\mu= 0.65$, the predicted values of $\lambda$ by the spline and the ANN are $0.2540$ and $0.0909$, respectively. The $L^{\infty}$ and $L^2$ error plots of the M\"untz spectral method using parameter values $\lambda=1$, $0.2540$ and $0.0909$ are shown in Figs. \ref{Plot2_1}, \ref{Plot2_2} and \ref{Plot2_3}, respectively. 
Fig. \ref{Plot2_1} illustrates that when $\lambda=1$, the approach essentially reduces to the classical spectral method, which cannot capture singularities and thus exhibits low accuracy. In contrast, Fig. \ref{Plot2_2} shows that adopting the spline-predicted parameter $\lambda= 0.2540$ achieves algebraic convergence and delivers higher accuracy than the classical spectral method. Finally, as evidenced in Fig. \ref{Plot2_3}, using the ANN-predicted parameter $\lambda = 0.0909$ yields exponentially decaying errors with respect to the polynomial degree $N$, demonstrating a substantial improvement in accuracy over previous approaches.

In summary, even in the relatively straightforward case of a one-to-one mapping from $\mu$ to $\lambda$, ANNs demonstrate certain advantages over spline interpolation in both efficiency and accuracy. Nevertheless, developing a simple ANN to predict optimal parameter values for the M\"ntz spectral method in time-fractional PDEs is not the central goal of our work. Our primary objective is to establish a fundamental optimization framework that incorporates ANNs into numerical solvers and trains them by minimizing solution errors. This framework is flexible and will be extended to more complex applications in the future, where $\lambda$ may be a vector $\lambda = (\lambda_1, \lambda_2, \ldots, \lambda_N)$ rather than a scalar, and where additional factors affecting solution's regularity can be introduced as inputs.
In such more complex settings, the mapping from the solution’s regularity to the optimal parameter(s) can become one-to-many, many-to-one, or even many-to-many. For instance, for a two-term time-fractional PDE, we can train an ANN to predict two parameters $(\lambda, q)$ for computing the exponents of the fractional polynomial basis functions as $\lambda_k = \lambda k + q$. This leads to a two-to-two mapping from $(\mu_1, \mu_2)$ to $(\lambda, q)$. In these more advanced applications, ANNs are expected to offer considerably greater benefits over traditional approaches such as spline or polynomial interpolation.



\section{Numerical results}
\label{sec4}

This section presents the numerical results of the machine learning optimized M\"untz spectral method for time-fractional PDEs in one and two dimensions, to demonstrate the effectiveness of the optimization.
Specifically, the accuracy improvement achieved through machine learning optimization is evaluated by comparing results from the M\"untz spectral method using FNN-predicted parameter values, the M\"untz spectral method with empirically selected parameter values, and the classical spectral method. 
The evaluation spans a variety of scenarios, including those with manufactured exact solutions involving both rational and irrational fractional powers, as well as a case where the exact solution is unknown.

For problems with smooth solutions, the classical spectral method can be directly applied by setting $\lambda = 1$ within the M\"untz spectral framework. However, since time-fractional PDEs typically yield solutions with low regularity, our numerical experiments are specifically designed to address cases involving these equations.




%

The $L^2$ and $L^\infty$ errors of the numerical solutions are defined as
\begin{equation*}
	\left\|u_e-u_L\right\|_{L^2}=\left(\dfrac{1}{K}\sum^K_{k=1}\left|u_e(x_k,t_{end})-u_L(x_k,t_{end})\right|^2\right)^{1/2},
\end{equation*}
and 
\begin{equation*}
	\left\|u_e-u_L\right\|_{L^\infty}=\max_k\left|u_e(x_k,t_{end})-u_L(x_k,t_{end})\right|,
\end{equation*}
where $t_{end}$ is the final time and $K$ is the number of nodes in space.  
The temporal and spatial errors are analyzed independently. To evaluate the temporal error, we use a sufficiently large value of $M$ to ensure that the spatial discretization error becomes negligible. Conversely, for the analysis of spatial error, we select a sufficiently large $N$ to minimize the impact of temporal discretization error.


\subsection{One-dimensional problems}
\begin{example}\label{exam1}
We consider the time-fractional convection-diffusion equation \eqref{TFDE} with a manufactured exact solution $u(x,t)=\sin(\pi x)t^{1-\mu}$. The time-fractional order is $\mu= 3/25$. 
\end{example}

\begin{figure}[htbp!]
\centering
\begin{subfigure}{0.48\textwidth}
	\includegraphics[width=\linewidth]{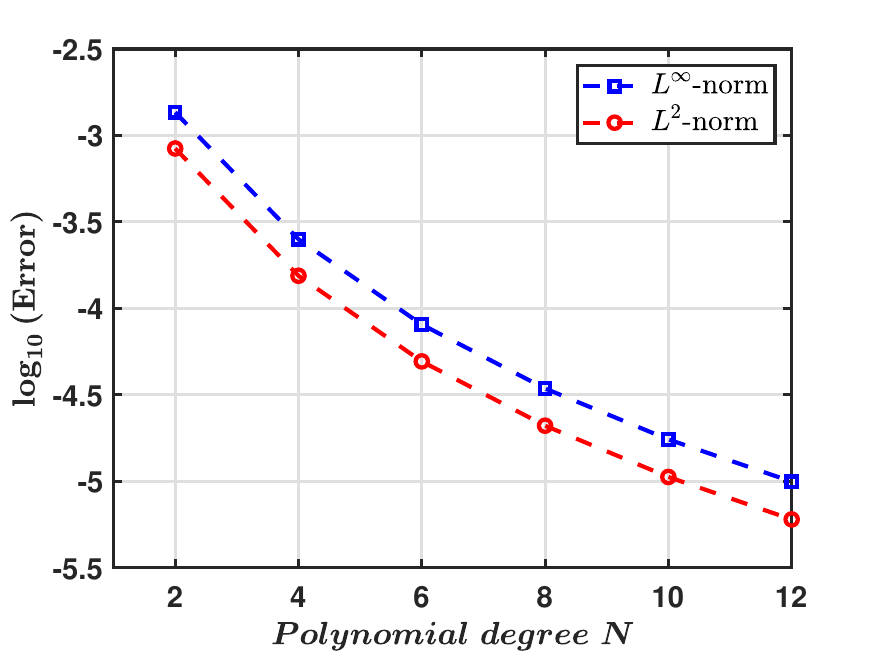}
	\caption{$\lambda=1$}
	\label{ErrPlot1t_1}
\end{subfigure}
\begin{subfigure}{0.48\textwidth}
	\includegraphics[width=\linewidth]{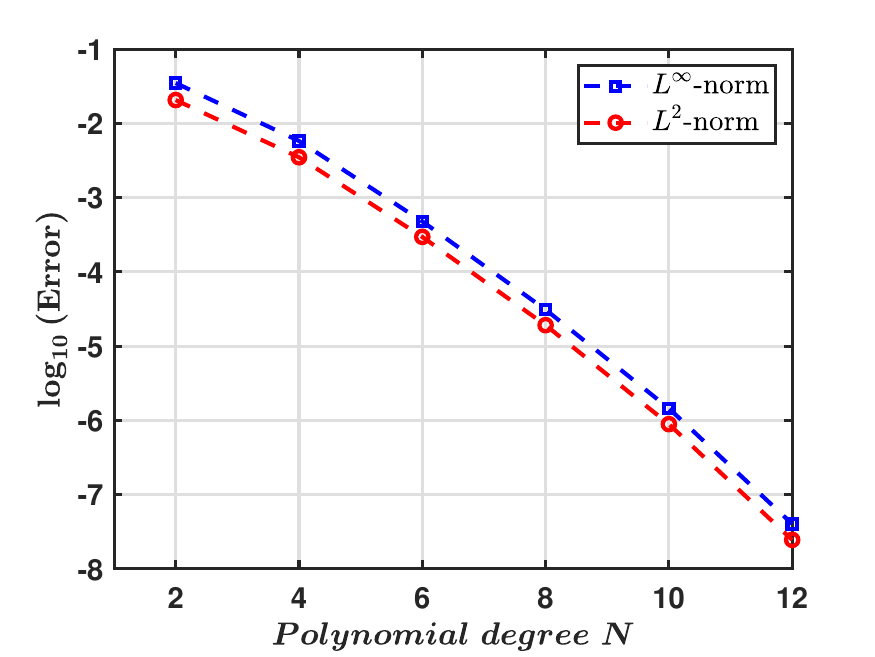}
	\caption{$\lambda=1/25$}
	\label{ErrPlot1t_2}
\end{subfigure}
\begin{subfigure}{0.48\textwidth}
	\includegraphics[width=\linewidth]{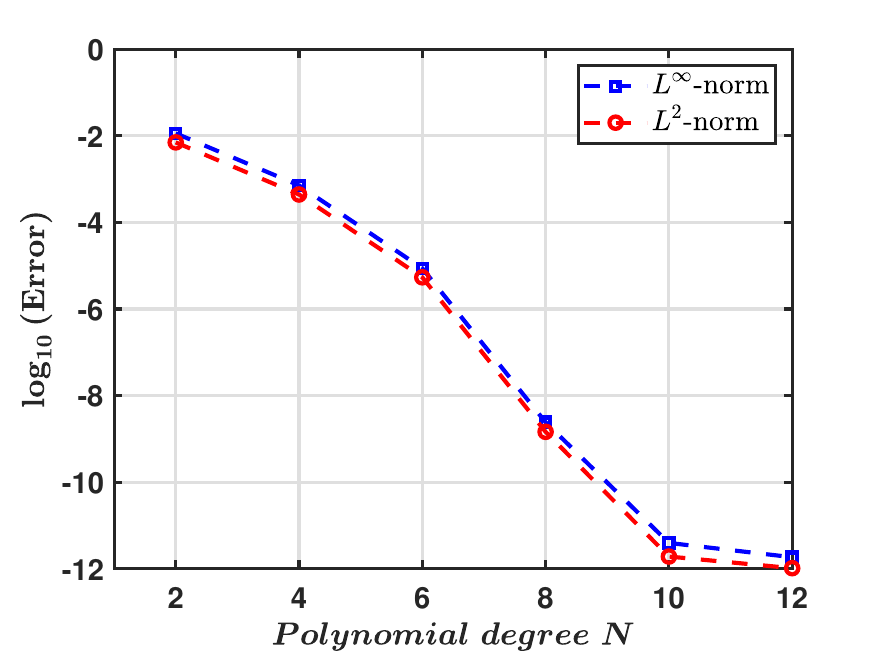}
	\caption{$\lambda=0.0962$}
	\label{ErrPlot1t_3}
\end{subfigure}
\caption{(Errors in time) $L^{\infty}$ and $L^2$ errors in semi-log scale versus $N$ with $\mu=3/25$ for different $\lambda$.}
\label{ErrPlot1t}
\end{figure}

We begin by analyzing the temporal discretization error with the number of spatial basis functions set to $M = 20$.
The $L^{\infty}$ and $L^2$ error plots of the M\"untz spectral method with $\lambda=1$, $1/25$ and $0.0962$ are shown in Figs. \ref{ErrPlot1t_1}, \ref{ErrPlot1t_2} and \ref{ErrPlot1t_3}, respectively. The plots in Fig. \ref{ErrPlot1t_1} show that, the M\"untz spectral method with $\lambda=1$ which is actually the classical spectral method, has very low accuracy since traditional polynomial basis functions are not capable of capturing singularities.
The plots in Fig. \ref{ErrPlot1t_2} are error curves of the M\"untz spectral method with a parameter value $\lambda=1/25$, selected based on the findings in \cite{Hou2018Muntz,Hou2017A}. The strategy of this selection is that, if the fractional order $\mu$ is in the form of a rational number $p/q$, the value $\lambda=1/q$ is suggested for the M\"untz spectral method. However, as shown in Fig. \ref{ErrPlot1t_2}, this parameter selection does not lead to high accuracy, which agrees with the analysis for the scenario when $q$ is large in Section \ref{issue}. 
The error plots of the M\"untz spectral method using the FNN predicted parameter value $\lambda= 0.0962$ are presented in Fig. \ref{ErrPlot1t_3}. It is shown that the errors decay exponentially with the polynomial degree $N$, demonstrating the significant accuracy improvement by using the FNN parameter prediction.  

\begin{figure}[htbp!]
\centering
\begin{subfigure}{0.48\textwidth}
	\includegraphics[width=\linewidth]{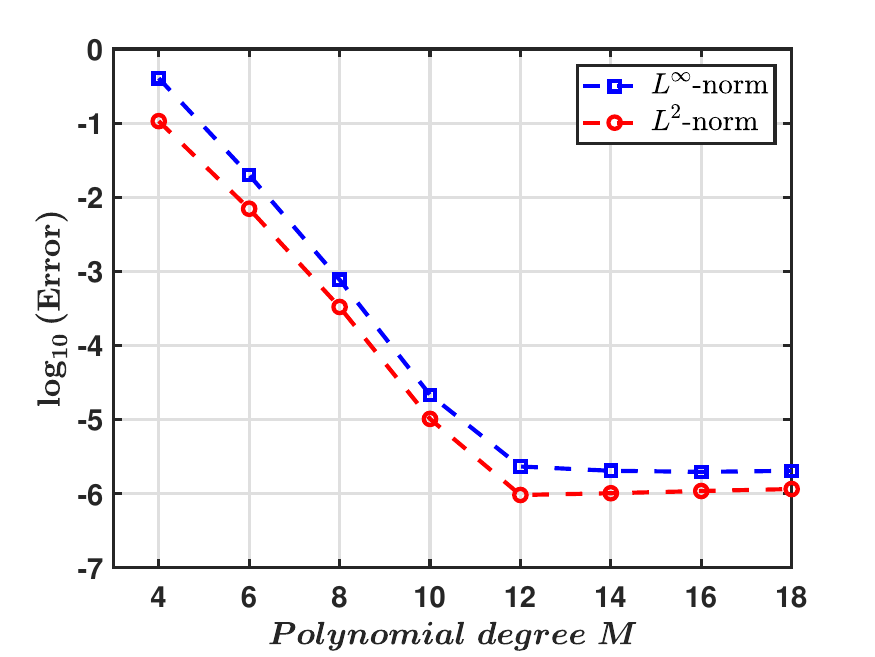}
	\caption{$\lambda=1$}
	\label{ErrPlot1s_1}
\end{subfigure}
\begin{subfigure}{0.48\textwidth}
	\includegraphics[width=\linewidth]{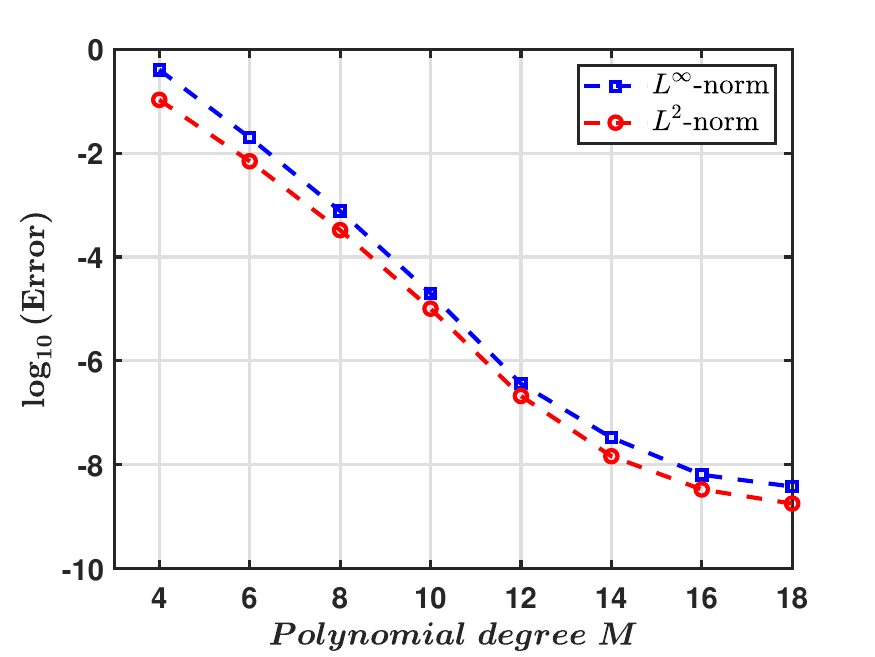}
	\caption{$\lambda=1/25$}
	\label{ErrPlot1s_2}
\end{subfigure}
\begin{subfigure}{0.48\textwidth}
	\includegraphics[width=\linewidth]{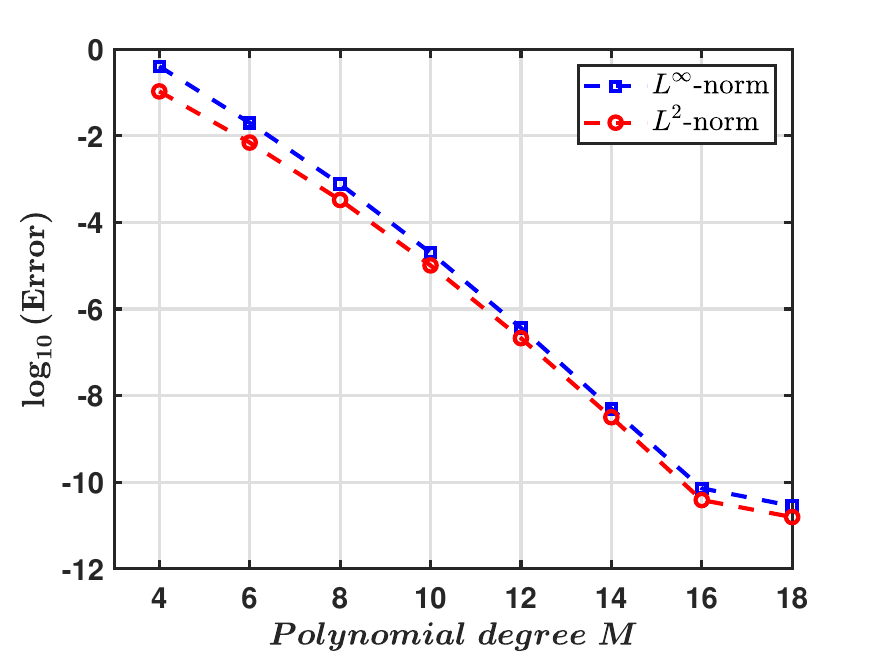}
	\caption{$\lambda=0.0962$}
	\label{ErrPlot1s_3}
\end{subfigure}
\caption{(Errors in space) $L^{\infty}$ and $L^2$ errors in semi-log scale versus $M$ with $\mu=3/25$ for different $\lambda$.}
\label{ErrPlot1s}
\end{figure}

We perform an analysis of the spatial discretization error with the number of temporal basis functions set to $N = 20$.
Figs. \ref{ErrPlot1s_1}, \ref{ErrPlot1s_2} and \ref{ErrPlot1s_3} plot $L^{\infty}$ and $L^2$ errors of the M\"untz spectral method with $\lambda=1$, $1/25$ and $0.0962$, respectively. Similar conclusions can be drawn from the results presented in Fig. \ref{ErrPlot1s} as those in Fig. \ref{ErrPlot1t}, suggesting that, in the space-time spectral method, the numerical accuracy in the time direction affects the numerical accuracy in the spatial direction.
Accordingly, in the remaining numerical experiments in one dimension, we will focus on presenting error results related to the temporal dimension, while omitting spatial error results. This approach is justified by the singular behavior exhibited by the solutions of time-fractional PDEs in the temporal dimension.

\begin{example}\label{exam2}
We consider the time-fractional convection-diffusion equation \eqref{TFDE} with a manufactured exact solution $u(x,t)=\sin(\pi x)t^{1+\mu}$. The time-fractional order is an irrational number $\mu= \sqrt{2}/2$. 
An adequate number of spatial basis functions is ensured by setting $M = 20$.
\end{example}

\begin{figure}[htbp!]
\centering
\begin{subfigure}{0.48\textwidth}
	\includegraphics[width=\linewidth]{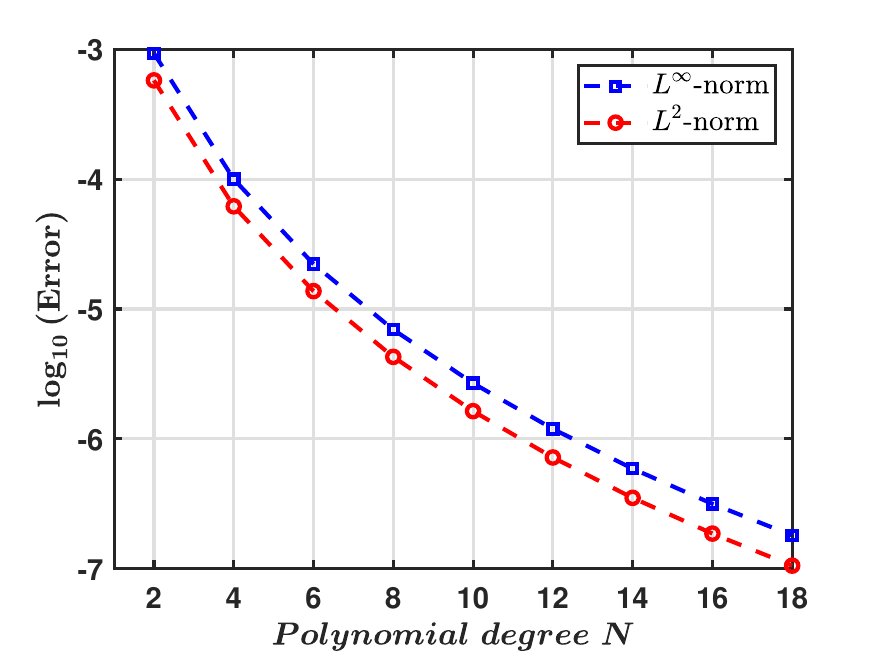}
	\caption{$\lambda=1$}
	\label{ErrPlot2_1}
\end{subfigure}
\begin{subfigure}{0.48\textwidth}
	\includegraphics[width=\linewidth]{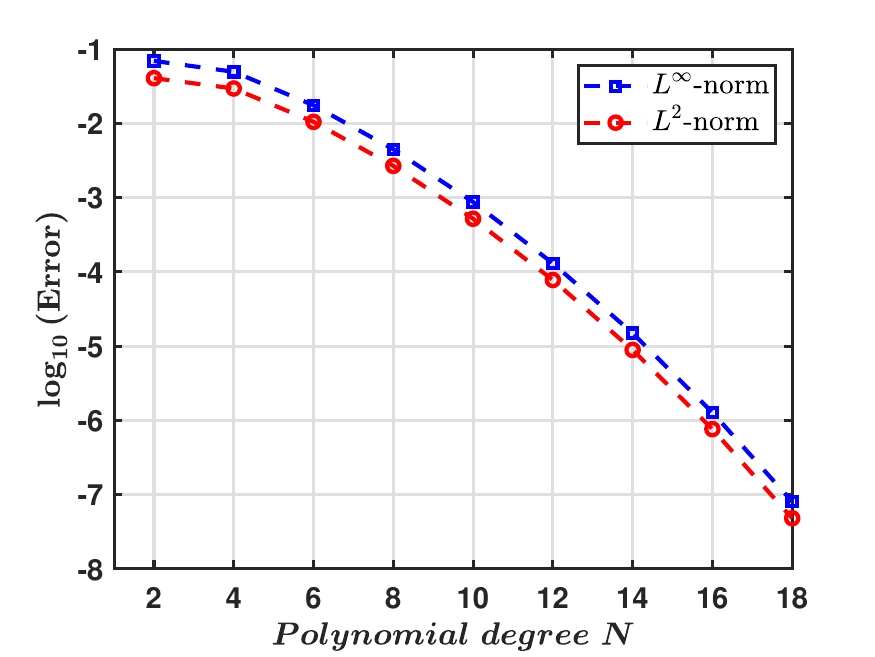}
	\caption{$\lambda=1/20$}
	\label{ErrPlot2_2}
\end{subfigure}
\begin{subfigure}{0.48\textwidth}
	\includegraphics[width=\linewidth]{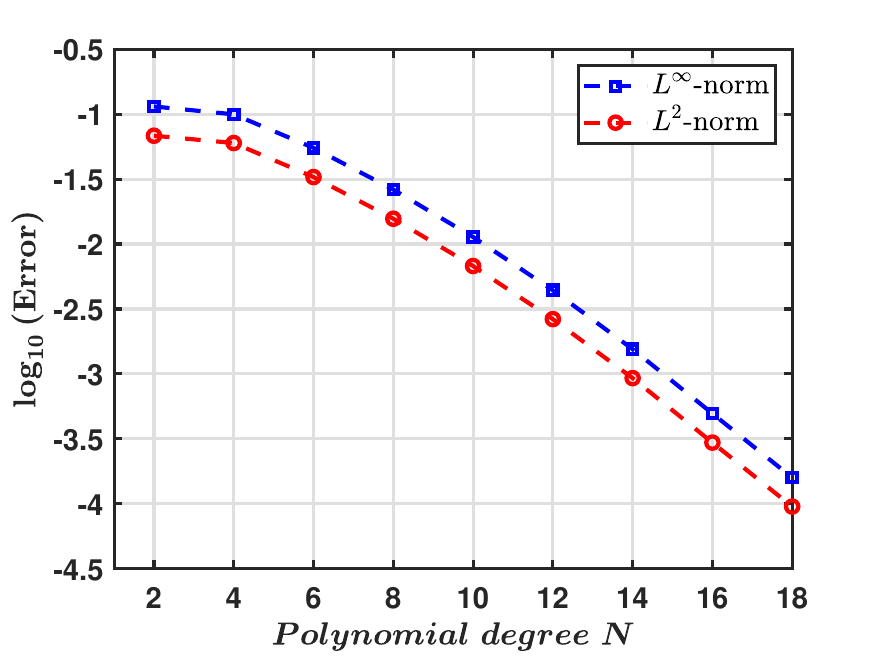}
	\caption{$\lambda=1/50$}
	\label{ErrPlot2_3}
\end{subfigure}
\begin{subfigure}{0.48\textwidth}
	\includegraphics[width=\linewidth]{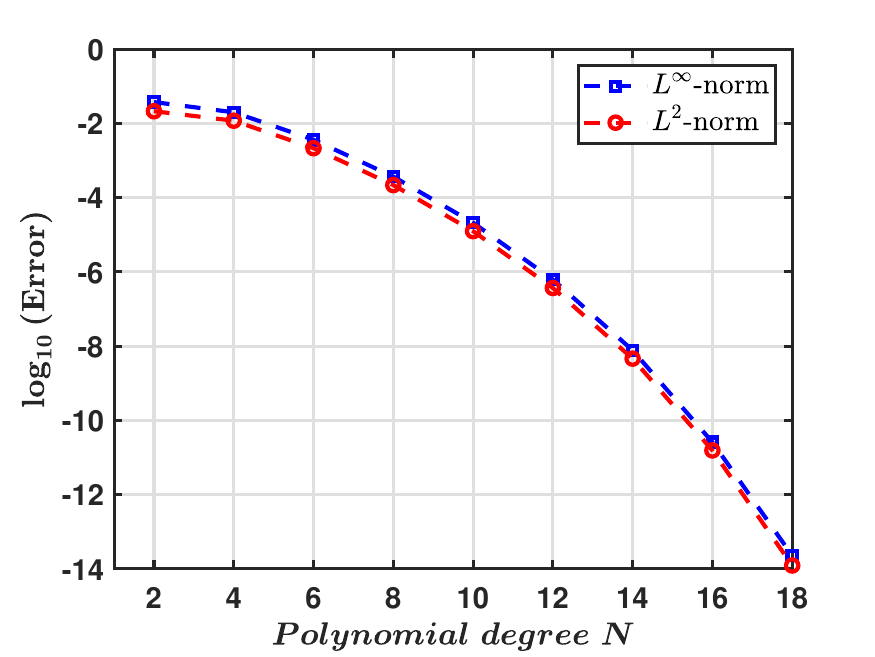}
	\caption{$\lambda=0.0899$}
	\label{ErrPlot2_4}
\end{subfigure}
\caption{(Errors in time) $L^{\infty}$ and $L^2$ errors in semi-log scale versus $N$ with $\mu=\sqrt{2}/2$ for different $\lambda$.}
\label{ErrPlot2}
\end{figure}

Error plots for the M\"untz spectral method with $\lambda=1$, $1/20$, $1/50$ and $0.0899$ are shown in Figs. \ref{ErrPlot2_1}, \ref{ErrPlot2_2}, \ref{ErrPlot2_3} and \ref{ErrPlot2_4}, respectively. 
As illustrated in Fig. \ref{ErrPlot2_1}, the classical spectral method ($\lambda = 1$) exhibits very low accuracy.
Given that the time-fractional order $\mu$ is irrational, it is challenging to select a $\lambda$ that yields a sufficiently smooth mapped solution $u(x,t^{1/\lambda})$. A common empirical approach is to choose $\lambda = 1/q$ with reasonably large integers $q$, such as 20 and 50, as tested here. However, as shown in Figs. \ref{ErrPlot2_2} and \ref{ErrPlot2_3}, these values do not achieve accurate results. Notably, the results in Fig. \ref{ErrPlot2_3} are even less accurate than those in Fig. \ref{ErrPlot2_1}, underscoring that $\lambda$ should not be excessively small, as discussed in Section \ref{issue}.
The error plots for the M\"untz spectral method using the FNN predicted $\lambda= 0.0899$ are presented in Fig. \ref{ErrPlot2_4}, showing exponential error decay with increasing polynomial degree $N$. This demonstrates a substantial improvement in accuracy with the FNN parameter prediction.

\begin{example}\label{exam3}
In this experiment, we evaluate the machine learning-optimized M\"untz spectral method by applying it to a general time-fractional convection-diffusion equation \eqref{TFDE} with a specified smooth forcing function. This represents a more complex and general case with an unknown exact solution but is expected to exhibit a singularity at $t=0$.
The forcing function is defined as $f(x,t)=\sin(\pi x)\sin(\pi t)$.
To ensure an adequate number of spatial basis functions, we set $M = 40$.
\end{example}

%

\begin{figure}[htbp!]
\centering
\begin{subfigure}{0.48\textwidth}
	\includegraphics[width=\linewidth]{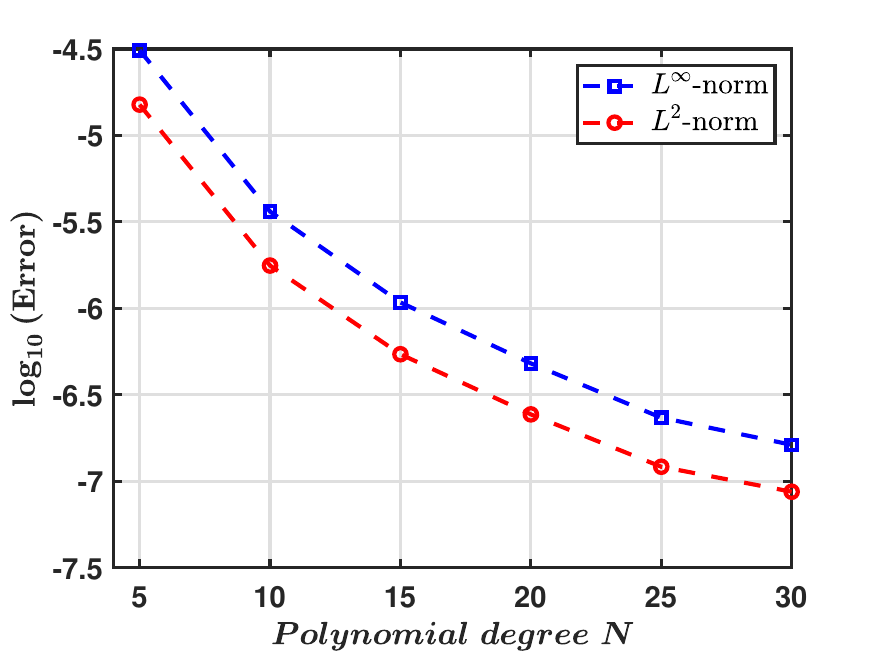}
	\caption{$\mu=0.3$, $\lambda=1$}
	\label{ErrPlot3_3}
\end{subfigure}
\begin{subfigure}{0.48\textwidth}
	\includegraphics[width=\linewidth]{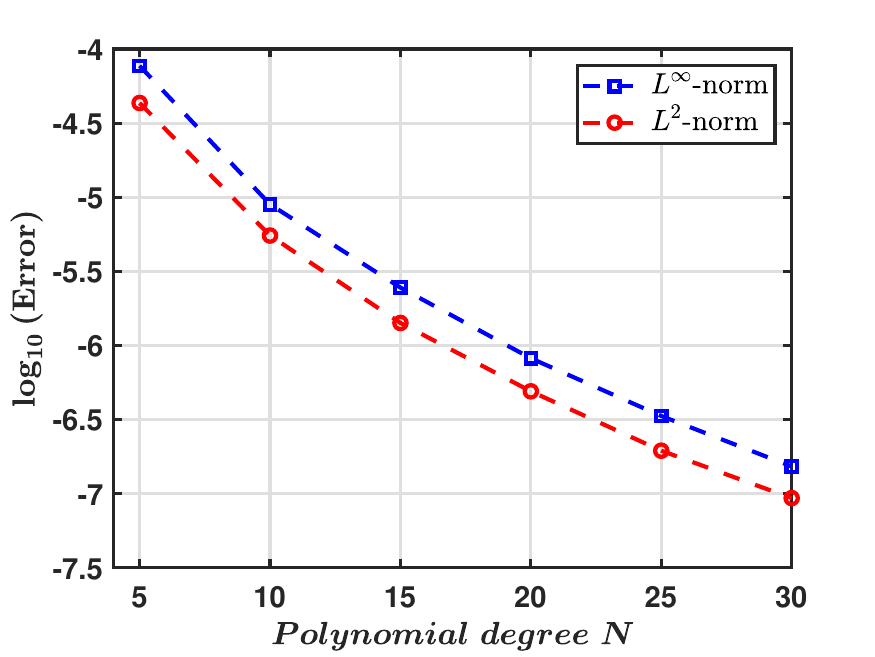}
	\caption{$\mu=0.5513$, $\lambda=1$}
	\label{ErrPlot3_4}
\end{subfigure}
\begin{subfigure}{0.48\textwidth}
	\includegraphics[width=\linewidth]{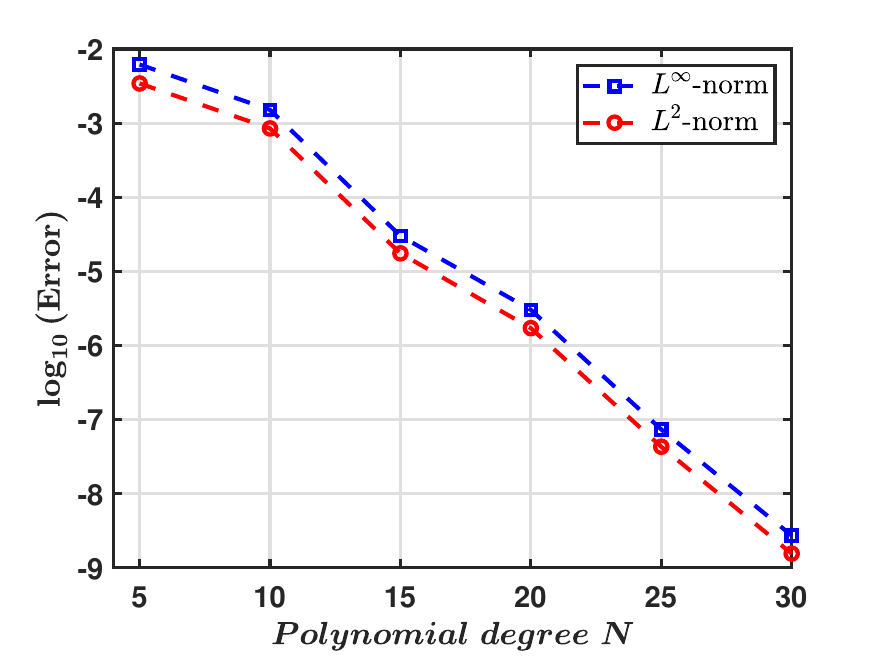}
	\caption{$\mu=0.3$, $\lambda=0.0951$}
	\label{ErrPlot3_1}
\end{subfigure}
\begin{subfigure}{0.48\textwidth}
	\includegraphics[width=\linewidth]{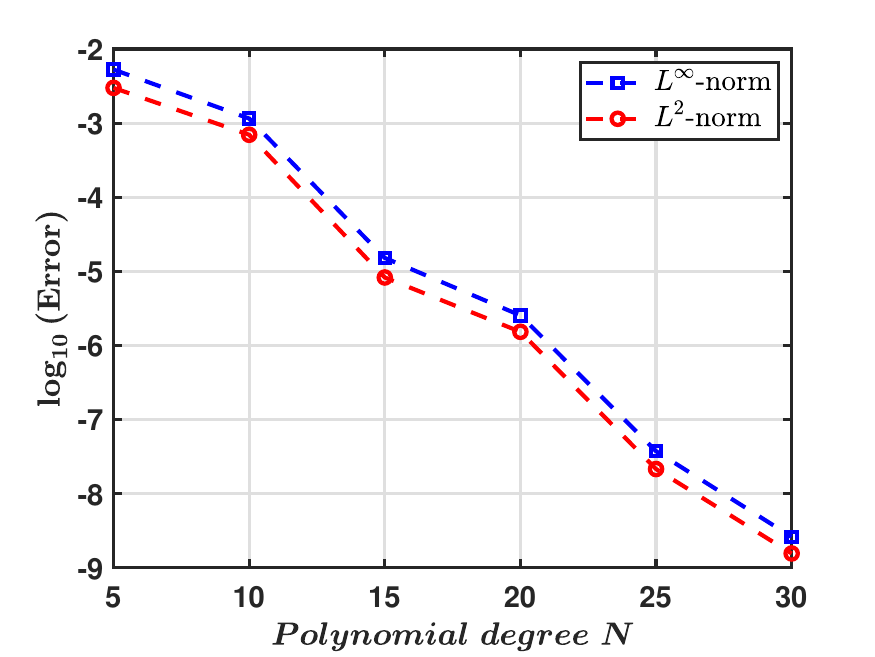}
	\caption{$\mu=0.5513$, $\lambda=0.0925$}
	\label{ErrPlot3_2}
\end{subfigure}
\caption{(Errors in time) $L^{\infty}$ and $L^2$ errors in semi-log scale versus $N$.}
\label{ErrPlot3}
\end{figure}

Fig. \ref{ErrPlot3} presents error plots for the M\"untz spectral method applied to two distinct time-fractional orders, $\mu=0.3$ and $0.5513$. 
The second time-fractional order is randomly selected within the interval $(0,1)$. 
For the two time-fractional orders, the parameter values predicted by the FNN are $\lambda = 0.0951$ and $0.0925$, respectively.
The numerical solutions computed using the machine learning optimized M\"untz spectral method with sufficient number of basis functions in both space and time by setting $M=N=40$, are used as the \textquote{exact} solutions to compute errors of the numerical solutions.
Figs. \ref{ErrPlot3_3} and \ref{ErrPlot3_4} illustrate the limited algebraic convergence when $\lambda = 1$, highlighting that solutions to time-fractional PDEs exhibit low regularity with respect to the time variable, which constrains the accuracy of the classical spectral method. By contrast, as shown in Figs. \ref{ErrPlot3_1} and \ref{ErrPlot3_2}, the machine learning optimized M\"untz spectral method, while still achieving algebraic convergence, delivers accuracy that is two orders of magnitude higher than that of the classical spectral method. This underscores its robustness and improved performance.

It is worth noting that in this general case, it is difficult to achieve spectral accuracy, posing a significant challenge. Ongoing research is being conducted to address this challenge, and we encourage readers to stay tuned for our future work on this topic.


\subsection{Two-dimensional problem}

\begin{example}\label{exam4}
In this test, we assess the performance of the machine learning optimized M\"untz spectral method on a two-dimensional sub-diffusion problem. This problem is modeled by the time-fractional convection-diffusion equation \eqref{TFDE} with a convection coefficient set to $\rho = 0$. The time-fractional order is chosen as $\mu = \sqrt{2}/2$, for which the FNN predicts a corresponding parameter value of $\lambda = 0.0899$. The manufactured exact solution is defined as  
\begin{equation*}
	u(x, y, t) = t^\eta \sin(\pi x) \sin(\pi y).
\end{equation*}
We conduct tests on three cases with different values of $\eta$, specifically, Case i: $\eta = 1 + \mu$, Case ii: $\eta = 1 - \mu$, and Case iii: $\eta = 3/7$. The time interval in this test is $I=[0,2]$.

\end{example}

\begin{figure}[htbp!]
	\centering
	\begin{subfigure}{\textwidth}
		\begin{minipage}{\linewidth}
			\includegraphics[width=0.45\linewidth]{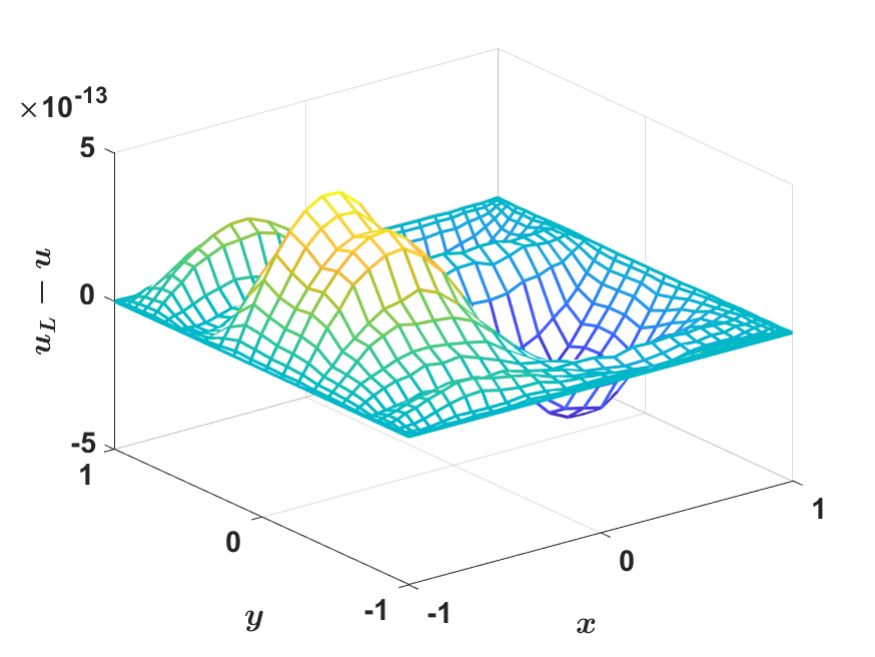}
			\includegraphics[width=0.45\linewidth]{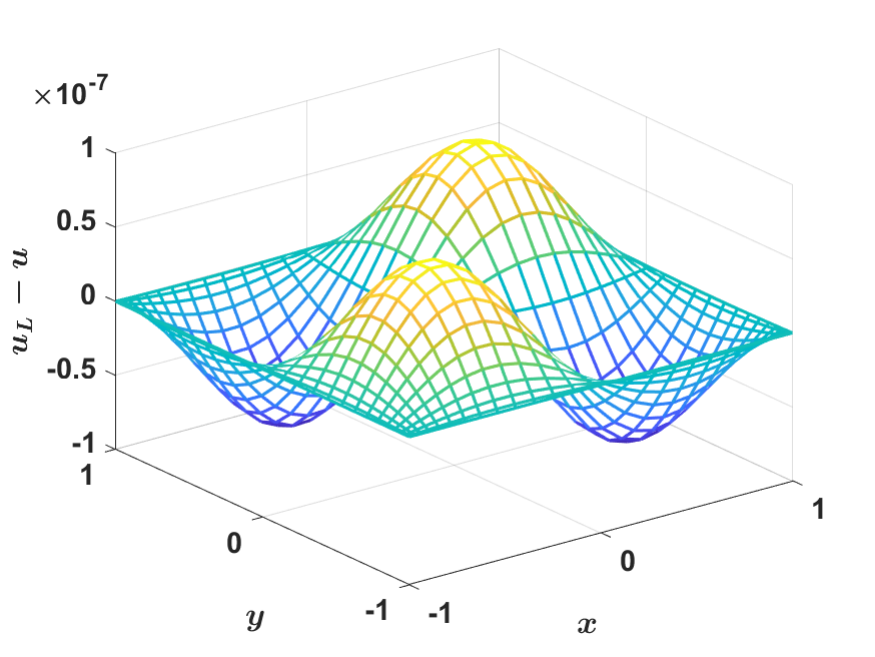}
		\end{minipage}\\
		\caption{Case i: $\lambda = 0.0899$ (left); $\lambda = 1$ (right)}
		\label{Plot2d4a_1} 
	\end{subfigure}
	
	\begin{subfigure}{\textwidth}
		\begin{minipage}{\linewidth}
			\includegraphics[width=0.45\linewidth]{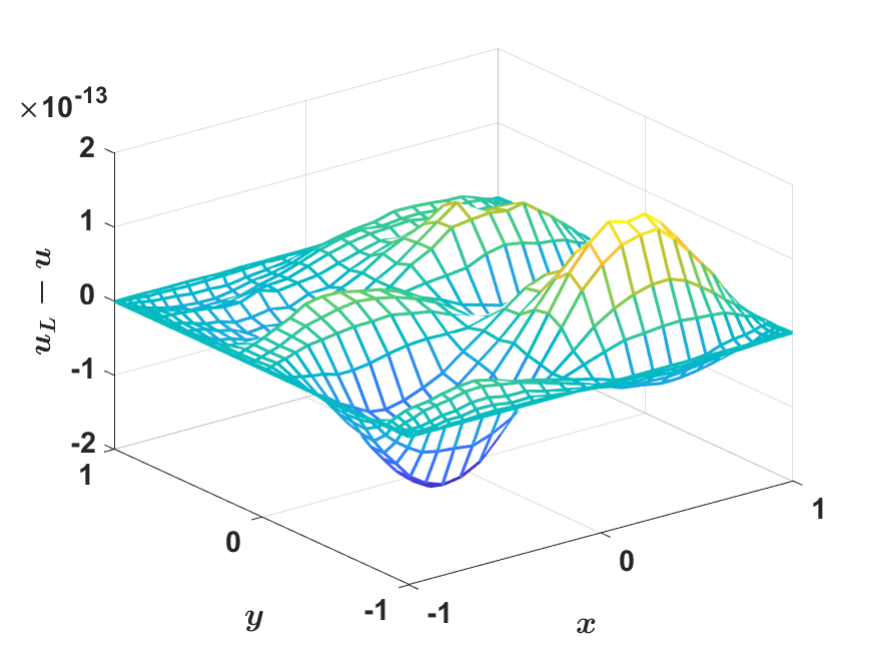}
			\includegraphics[width=0.45\linewidth]{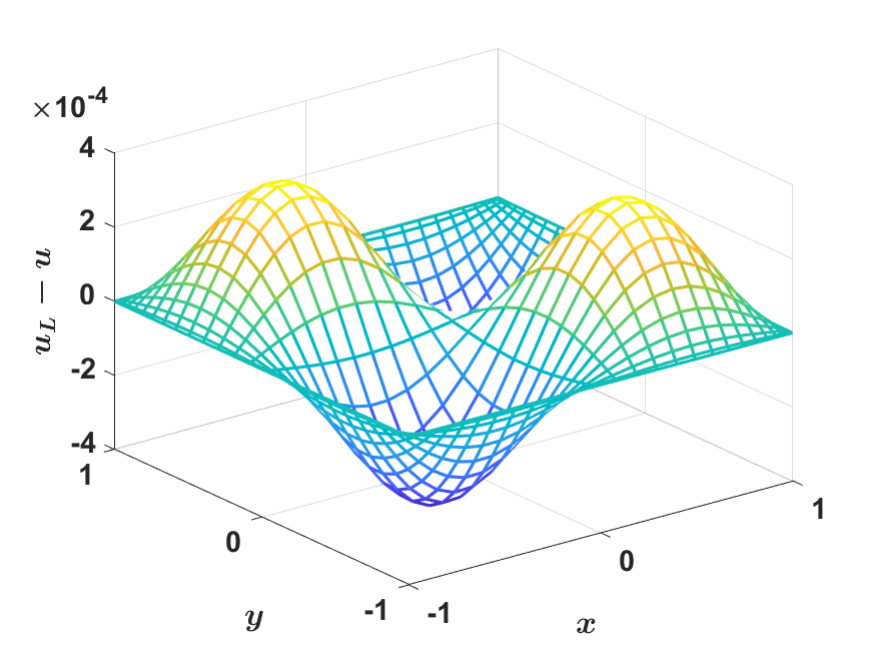}
		\end{minipage}\\
		\caption{Case ii: $\lambda = 0.0899$ (left); $\lambda = 1$ (right)}
		\label{Plot2d4a_2} 
	\end{subfigure}
	
	\begin{subfigure}{\textwidth}
		\begin{minipage}{\linewidth}
			\includegraphics[width=0.45\linewidth]{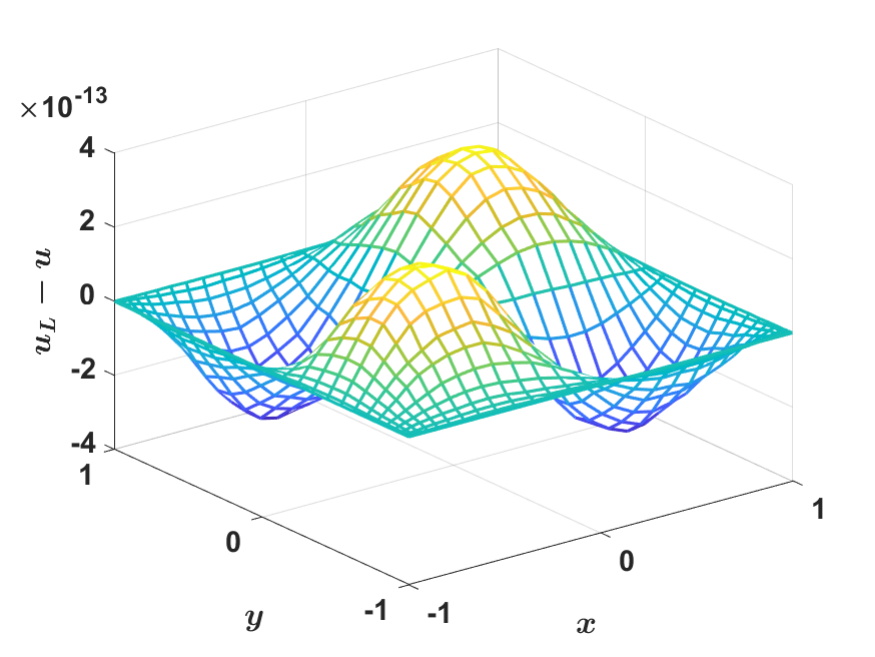}
			\includegraphics[width=0.45\linewidth]{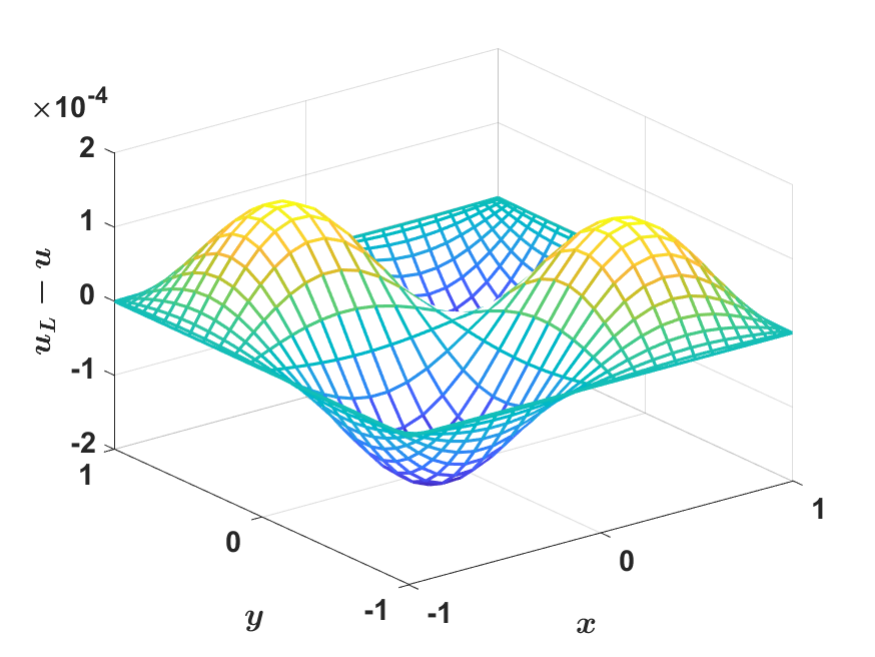}
		\end{minipage}\\
		\caption{Case iii: $\lambda = 0.0899$ (left); $\lambda = 1$ (right)}
		\label{Plot2d4a_3} 
	\end{subfigure}
	
	\caption{Case i to Case iii: Error $u_L-u$ at $t=2$.}
	\label{Plot2d_4a} 
\end{figure}

\begin{figure}[htbp!]
\centering
\begin{subfigure}{0.48\textwidth}
	\includegraphics[width=\linewidth]{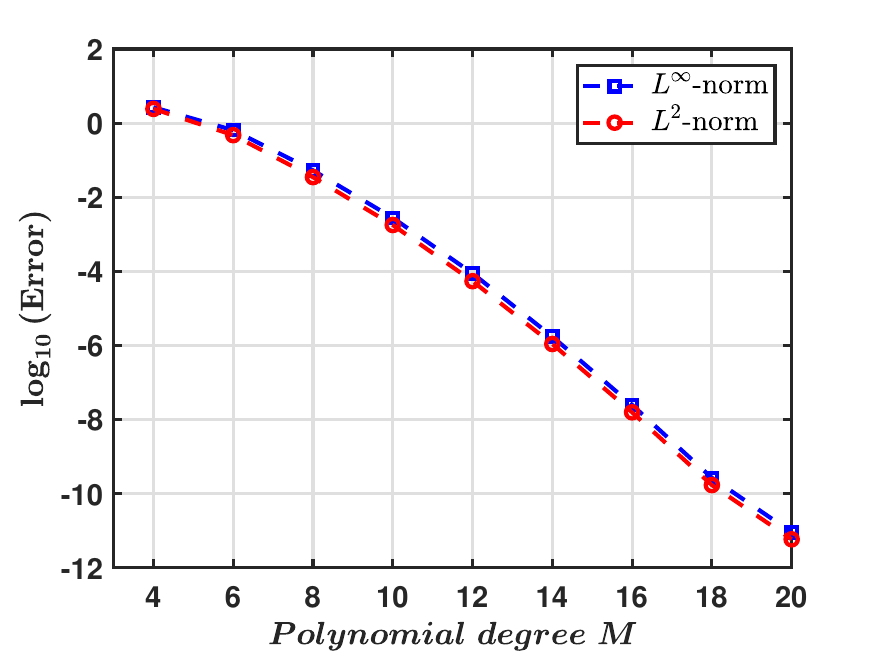}
	\caption{(Errors in space) $\eta=1+\mu$}
	\label{Plot2d4_1}
\end{subfigure}
\begin{subfigure}{0.48\textwidth}
	\includegraphics[width=\linewidth]{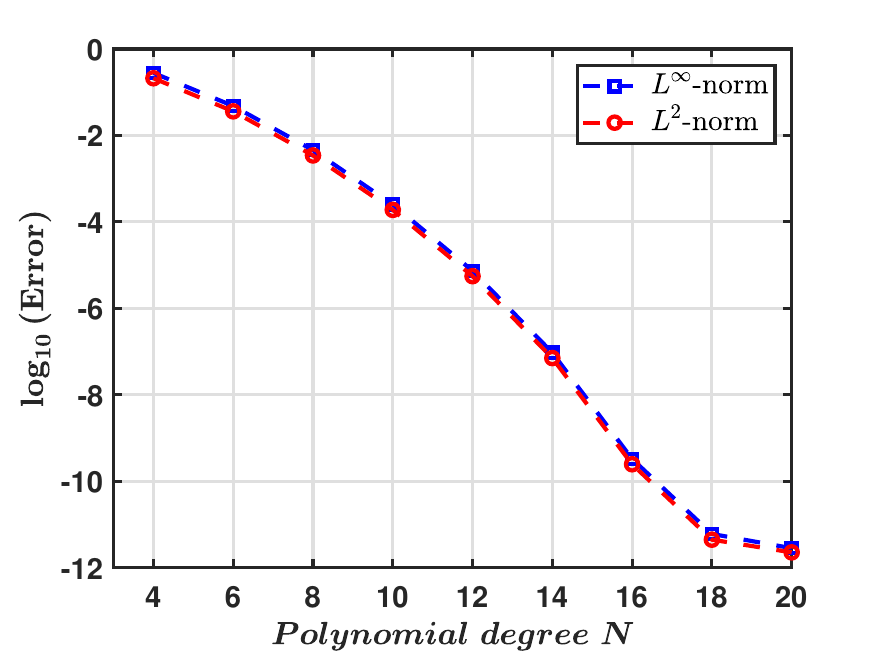}
	\caption{(Errors in time) $\eta=1+\mu$}
	\label{Plot2d4_2}
\end{subfigure}
\begin{subfigure}{0.48\textwidth}
	\includegraphics[width=\linewidth]{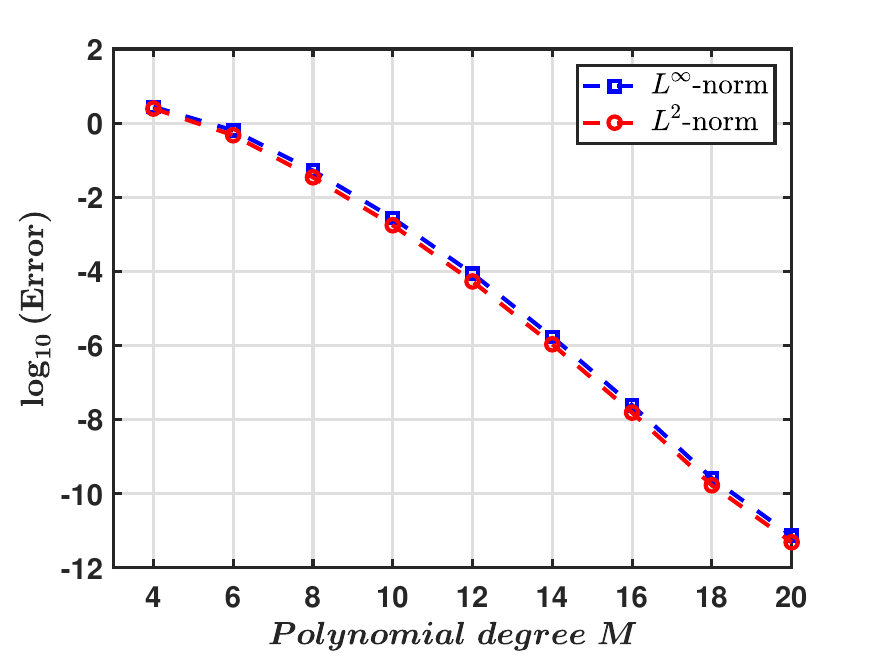}
	\caption{(Errors in space) $\eta=1-\mu$}
	\label{Plot2d4_3}
\end{subfigure}
\begin{subfigure}{0.48\textwidth}
	\includegraphics[width=\linewidth]{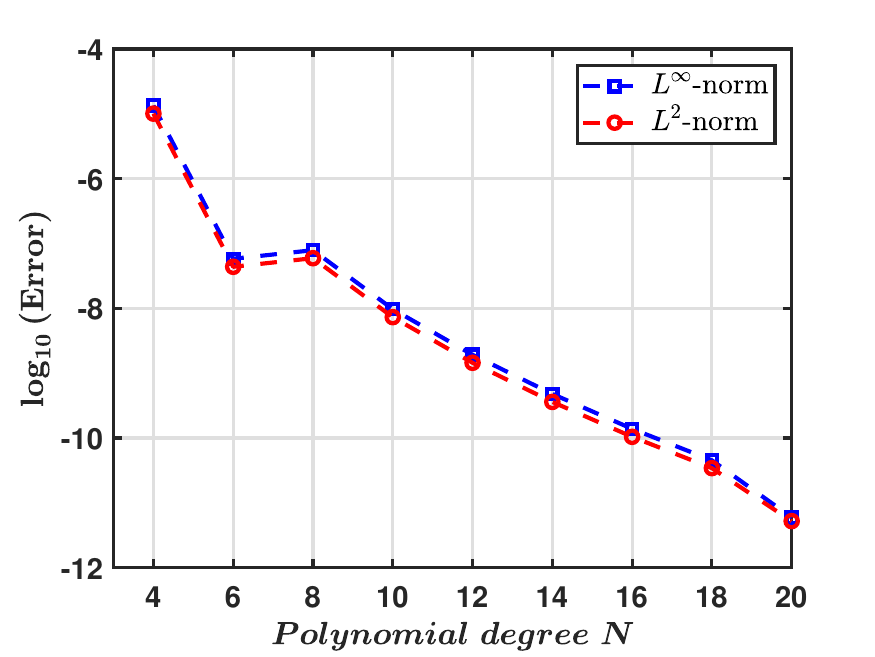}
	\caption{(Errors in time) $\eta=1-\mu$}
	\label{Plot2d4_4}
\end{subfigure}
\begin{subfigure}{0.48\textwidth}
	\includegraphics[width=\linewidth]{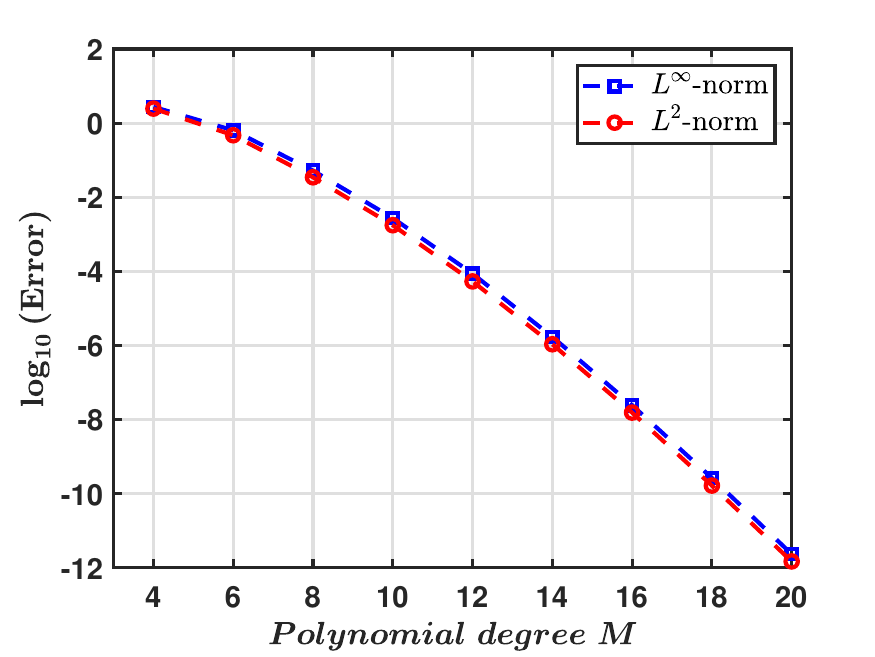}
	\caption{(Errors in space) $\eta=3/7$}
	\label{Plot2d4_5}
\end{subfigure}
\begin{subfigure}{0.48\textwidth}
	\includegraphics[width=\linewidth]{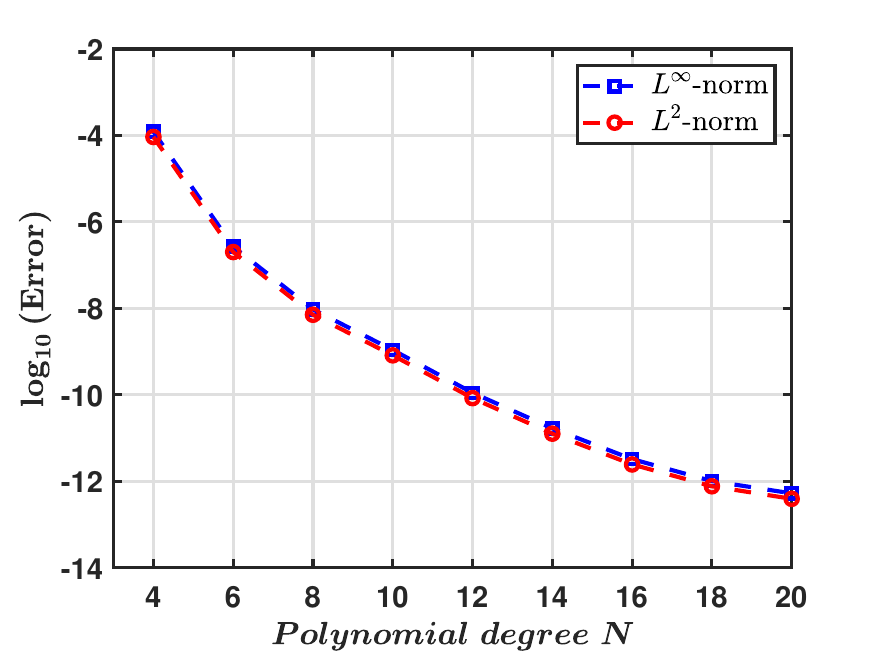}
	\caption{(Errors in time) $\eta=3/7$}
	\label{Plot2d4_6}
\end{subfigure}
\caption{Error plots in semi-log scale with $\mu=\sqrt{2}/2$ for different $\eta$.}
\label{Plot2d4}
\end{figure}

Figs. \ref{Plot2d4a_1} (left)-\ref{Plot2d4a_3} (left) illustrate the corresponding solution errors of the machine learning optimized M\"untz spectral method at the final time $t= 2$ for the three test cases. As shown, the pointwise errors are on the order of $10^{-13}$, underscoring the high accuracy of the numerical solutions.
For comparison, Figs. \ref{Plot2d4a_1} (right)-\ref{Plot2d4a_3} (right) show the solution errors of the classical spectral method, which exhibits limited accuracy across all three cases.
Additionally, Fig. \ref{Plot2d4} presents error plots of the machine learning optimized M\"untz spectral method, demonstrating spectral convergence of the solution error in both space and time.
For this analysis, $N=30$ and $M=30$ are used to assess spatial and temporal errors, respectively. 

The numerical results indicate that the FNN-based parameter prediction model, initially trained in one dimension, generalizes effectively to two-dimensional problems without requiring retraining. This generalization is expected, as time-fractional PDEs exhibit singularities only in the temporal domain, allowing the trained model to deliver consistent performance across varying spatial dimensions.

\section{Conclusions}
\label{sec5}

This paper introduces a machine learning optimization approach for the M\"untz spectral method, where a feedforward neural network (FNN) is trained to predict the optimal parameter $\lambda$ by minimizing solution errors of a one-dimensional time-fractional convection-diffusion equation. This approach effectively addresses the challenge of time-consuming parameter tuning. 
Numerical results demonstrate that the machine learning optimization significantly enhances the accuracy of the M\"untz spectral method. Furthermore, the trained network generalizes effectively to two-dimensional time-fractional problems, highlighting its robust performance across spatial dimensions.

While this paper centers on spectral methods for time-fractional PDEs, the proposed machine learning optimization framework is versatile and can be extended to numerical methods employing fractional polynomial bases for a broader range of differential equations. This approach holds promise for accurately solving problems where solutions exhibit singularities.
Future work will focus on expanding the proposed framework to accommodate more application scenarios and exploring more efficient neural network architectures to further enhance the accuracy and efficiency of parameter selection.



\bibliographystyle{siamplain}
\bibliography{references}
\end{document}